\pdfoutput=1
\documentclass[conference]{IEEEtran}

\IEEEoverridecommandlockouts

\usepackage{cite}
\usepackage{amsmath,amssymb,amsfonts}
\usepackage{algorithmic}
\usepackage{graphicx}
\usepackage{booktabs}
\usepackage{xcolor}
\usepackage{pgfplots}
\pgfplotsset{compat=1.16}
\usepackage{url}
\usepackage[hidelinks]{hyperref}

\graphicspath{{figures/}}

\newcommand{\GFtwo}{\ensuremath{\mathrm{GF}(2)}}
\newcommand{\LN}[2]{\ensuremath{L_N[#1,#2]}}

\newcommand{\AIStatement}{Generative AI and agentic coding tools were used.
Software development: Anthropic Claude Haiku~4.5,
Sonnet~4.5--5.5, Opus~4.6--5.5 and Fable~5--5.1 (via Claude Code); OpenAI GPT-5
and GPT-5.2; Google Gemini~2.5~Pro, 3~Pro and 3.1~Pro. Data and manuscript preparation: Claude Sonnet~5, Opus~5 and 5.5, and Fable~5 and 5.1 (via
Claude Code). All AI-generated output was manually reviewed and verified by the authors, who retain ownership and take full responsibility for the software and the content of this paper.}

  \newcommand{\hpcsite}{the PC${}^2$}
  \newcommand{\swname}{CUDA-MPQS}
  \newcommand{\papertitle}{\swname: A GPU-Resident Self-Initializing\\
    Quadratic Sieve, and the Factorization of RSA-155}
  \newcommand{\pdfpapertitle}{CUDA-MPQS: A GPU-Resident Self-Initializing
    Quadratic Sieve, and the Factorization of RSA-155}
  \newcommand{\pdfpaperauthor}{Fabian Januszewski; Christoph Heinrichs}
  \newcommand{\heinrichscite}{heinrichs2026}
  \newcommand{\geiselmanncite}{geiselmann2006}
\hypersetup{pdftitle={\pdfpapertitle}, pdfauthor={\pdfpaperauthor},
  pdfsubject={}, pdfkeywords={}}

\begin{document}

\title{\papertitle}

  \author{
  \IEEEauthorblockN{Fabian Januszewski}
  \IEEEauthorblockA{\textit{Institut f\"ur Mathematik} \\
  \textit{Paderborn University}\\
  Paderborn, Germany \\
  fabian.januszewski@math.uni-paderborn.de}
  \and
  \IEEEauthorblockN{Christoph Heinrichs}
  \IEEEauthorblockA{\textit{Institut f\"ur Mathematik} \\
  \textit{Paderborn University}\\
  Paderborn, Germany \\
  christoph.heinrichs@math.uni-paderborn.de}
  }

\maketitle

\begin{abstract}
The quadratic sieve (QS) is an irregular, cache-hostile, branch-heavy integer
factorization algorithm; prior GPU work accelerates individual stages of it. We present
\swname{}, a self-initializing quadratic sieve in which \emph{every} stage ---
polynomial initialization, sieving, relation post-processing, large-prime matching,
\GFtwo{} matrix construction, Block Wiedemann linear algebra and the square root ---
executes on the GPU, the CPU restricted to orchestration, setup and I/O, the sieve window
99.9\% GPU-busy with no device-wide or event synchronization and no synchronous copy.
With it we factored the 512-bit RSA-155 challenge modulus end to end on GPUs: a
64$\times$H100 sieve feeding a single-H100 \GFtwo{} solve, 700.6~GPU-h and 242~kWh in
total, 98.4\% of it sieving. To our knowledge this is the largest integer factored by
the quadratic sieve, exceeding the prior QS record RSA-150 (Buhrow, YAFU, 2025); all
comparable general factorizations since 1996, RSA-155's own in 1999 included, used the
number field sieve. On a single device we factor RSA-100 in 29.2~s on an H100 and
51.0~s on a consumer RTX~5070~Ti; in a controlled same-modulus measurement, with energy
measured on both sides, the H100 is faster than the fastest CPU quadratic sieve, YAFU, on 96 cores of an EPYC~9655 (106--122~s), and
than CADO-NFS's SIQS (292~s) and GNFS (263~s). Across four discrete GPUs (a
3.7$\times$ range in CUDA cores) RSA-100 throughput is close to proportional to core
count, memory bandwidth being a weaker predictor, and measured cost over seven
benchmark factorizations from RSA-100 to RSA-155 grows as $\LN{1/2}{1}^{1.04}$ ($R^2=0.993$). On
RSA-150, the identical modulus on which the prior QS record was set, we measure
302.9~GPU-h against that run's 11{,}664 core-hours.
\end{abstract}


\section{Introduction}
\label{sec:intro}

The quadratic sieve~\cite{pomerance1985qs,pomerance1996tale} is, from a
parallel-computing standpoint, an unpromising GPU workload: its inner loop scatters
byte-granular read-modify-write updates across an interval far larger than any cache,
over millions of primes each with its own stride; the smoothness test is branch-heavy
and data-dependent; and the linear algebra operates over \GFtwo{}, a field no vendor
library covers.
Accordingly, prior GPU work on factorization has accelerated \emph{stages} --- sieving,
ECM, NFS cofactorization, polynomial selection and the \GFtwo{} sparse matrix--vector
core~\cite{wang2011gpuqs,bernstein2009ecm,miele2014cofactorization,msieve,schmidt2011spmv,schmidt2013bwgpu,msievelacuda}
(\S\ref{sec:related}) --- while the host retains the rest, and with it the pipeline's
synchronization structure.
Whether the \emph{whole} pipeline can run efficiently on a GPU, and be competitive
with the best CPU quadratic sieves, has to our knowledge not been demonstrated.

In \swname{}, a self-initializing quadratic sieve (SIQS), polynomial initialization,
sieving, trial division, relation deduplication, large-prime matching, \GFtwo{} matrix
construction, Block Wiedemann and the square root all execute in device kernels. The
engineering consequence is a sieve loop that on a single GPU reaches a
zero-synchronization steady state: no blocking host sync in the loop body, and even
the decision to \emph{stop} made by a device kernel writing a structure the host polls
without synchronizing. The scientific consequence is a lopsided cost structure:
linear algebra, the stage hardest to keep resident, is 1.6\% of our largest run.
We measured the competitive consequence: at RSA-100 one H100 (v1.0.8) is $3.6$--$4.2\times$
faster than the fastest CPU quadratic sieve, YAFU on 96 cores of an EPYC~9655, and
$9.0$--$10.0\times$ faster than CADO-NFS's SIQS and GNFS, on the same modulus
(\S\ref{sec:cpubaseline}) --- to our knowledge the first GPU quadratic sieve whose
end-to-end wall-clock matches, and here exceeds, state-of-the-art CPU QS.

The quadratic sieve is asymptotically inferior to the number field sieve and the wrong
tool for record-setting general factorizations since 1996 (RSA-155 fell to GNFS in
1999~\cite{cavallar2000rsa155}; the record is now RSA-896 (2026)~\cite{weis2026rsa896});
nothing here narrows a security margin and we make no cryptanalytic claim. What the QS offers is a complete and unusually demanding
pipeline --- irregular sieving, hash-join-style relation matching, \GFtwo{} sparse
linear algebra --- on which ``can this class of algorithm be made GPU-resident, and does it scale?'' admits
a measured answer.

\subsection*{Contributions}

We contribute, first, \textbf{a fully GPU-resident QS pipeline with a
zero-synchronization sieve steady state} --- to our knowledge the first quadratic sieve
in which every stage executes on the GPU, the CPU restricted to orchestration, setup and I/O
(\S\ref{sec:design}, host residuals in \S\ref{sec:residency}) --- resting on a
\GFtwo{}-native auto-tuned SpMM engine of ten formats whose
dynamic-programming solver assigns \emph{different} formats to different density
segments of one matrix (\S\ref{sec:linalg}) and an on-device slab hash table for
large-prime matching (\S\ref{sec:lp}). Second, \textbf{the factorization of RSA-155
(512-bit) end to end on GPUs} (\S\ref{sec:rsa155}): 700.6~GPU-h, 242~kWh, sieve $=$
98.4\% of cost, to our knowledge the largest integer factored by the quadratic sieve
(\S\ref{sec:qsrecords}). Third, \textbf{a two-axis scaling characterization} --- seven
benchmark factorizations from RSA-100 to RSA-155 (\S\ref{sec:benchfact}) and five GPUs, the four
discrete ones spanning 3.7$\times$ in CUDA cores (\S\ref{sec:corescaling}) --- with a
same-modulus head-to-head at RSA-150 against the prior record holder
(\S\ref{sec:rsa150}), 302.9~GPU-h against 11{,}664 core-hours inside a cost band
committed to before the run, and a controlled CPU baseline at RSA-100, same
modulus, energy measured on both sides (\S\ref{sec:cpubaseline}).
 \section{Background}
\label{sec:bg}

We fix notation and state only what later sections depend on
(see~\cite{pomerance1985qs,silverman1987mpqs,contini1997siqs,crandall2005primes,\heinrichscite}).
To factor an odd composite $N$ it suffices to find $x^2 \equiv y^2 \pmod N$ with $x
\not\equiv \pm y$, whereupon $\gcd(N, x\pm y)$ splits $N$. The QS collects
\emph{relations}, values of a quadratic polynomial factoring completely over a
\emph{factor base} $\mathcal{F}$ of the primes $p<F$ with
$\left(\tfrac{N}{p}\right)=1$, as rows of exponent parities over \GFtwo{}; past
$|\mathcal{F}|$ of them any kernel vector gives a square on both
sides~\cite{dixon1981}.
\emph{Sieving}~\cite{davis1984factorization,pomerance1988pipeline,pomerance1996tale}
finds them: for each $p\in\mathcal{F}$ one adds $\log_2 p$ at every offset in the
progressions of the polynomial's roots mod $p$, and positions whose log-sum approaches
$\log_2|Q(x)|$ are candidates.

The multiple-polynomial variant~\cite{silverman1987mpqs} sieves the family
\begin{equation}
  Q_{(a,b)}(x) = (ax+b)^2 - N, \qquad b^2 \equiv N \!\!\pmod a ,
  \label{eq:poly}
\end{equation}
over $x\in[-M,M)$; with $b^2-N=ac$ one sieves $\widehat{Q}_{(a,b)}(x)=ax^2+2bx+c$,
re-attaches the known factorization of $a$, and takes $a\approx\sqrt{2N}/M$ to minimize
$\max_x|\widehat{Q}(x)|$. Self-initialization
(SIQS)~\cite{alford1995siqs,contini1997siqs,peralta1995hypercube} makes $a=q_1\cdots
q_k$ a product of $k$ factor-base primes (the \emph{hypercube dimension}); by CRT every
$b=\sum_i\pm B_i$ is a valid root, where
$B_i=(a/q_i)\,[\,t_i(a/q_i)^{-1}\bmod q_i\,]$ with
$t_i^2\equiv N \pmod{q_i}$; as $Q_{(a,-b)}(x)=Q_{(a,b)}(-x)$, this gives $2^{k-1}$
distinct polynomials per $a$, and Gray-code
enumeration~\cite{knuth2011taocp4a} makes each switch an $O(1)$ per-prime root update
($\pm 2a^{-1}B_i \bmod p$), which is the structure \S\ref{sec:sieve} exploits.

A cofactor below a bound $L>F$ may be left unfactored, yielding \emph{partials}; two
partials sharing a large prime multiply into a full relation~\cite{silverman1987mpqs},
and matches grow quadratically in the partial count, a birthday effect that drives
our design (\S\ref{sec:lp}) and complicates yield prediction in truncated autotuning
(\S\ref{sec:autotune}). We implement neither multiple-large-prime
variant~\cite{lenstra1994twolp,leyland2002threelp}.

The relation matrix is large and sparse (\S\ref{sec:rsa155}). We use Block
Wiedemann~\cite{coppersmith1994bw,wiedemann1986,kaltofen1995bw,villard1997bw} (Block
Lanczos~\cite{montgomery1995} being the alternative), whose three stages (Krylov generation, the generating polynomial or \emph{lingen},
reconstruction) are dominated by sparse \GFtwo{} matrix products (\S\ref{sec:linalg}).
 \section{Related Work}
\label{sec:related}

\subsection{GPU factorization, stage by stage}

Two entries bound our claim of \emph{end-to-end residency}. {QS sieving on GPUs
is not new}: Wang~\cite{wang2011gpuqs} reported GPU sieving experiments for the
quadratic sieve in 2011, citing an earlier GPU multiple-polynomial sieve with one
polynomial per thread.
{GPU \GFtwo{} block linear algebra is not new}: Schmidt, Aribowo and
Dang~\cite{schmidt2011spmv,schmidt2013bwgpu} accelerated Block Wiedemann's sparse
matrix--vector core on multi-GPU systems for NFS matrices (4--8$\times$ over
multicore, Berlekamp--Massey and orchestration staying on the host), and msieve's
\texttt{lacuda} line in the NFS@Home fork~\cite{msievelacuda} drives block
Lanczos~\cite{montgomery1995} with a CUDA \GFtwo{} SpMV engine, while mainline msieve
uses the GPU only for NFS polynomial selection~\cite{msieve}. What \S\ref{sec:linalg}
adds is ten storage formats, per-segment mixed assignment by a dynamic-programming
solver, and surrounding stages that stay on the device. Bernstein et
al.~\cite{bernstein2009ecm} put ECM on graphics cards and Miele et
al.~\cite{miele2014cofactorization} NFS cofactorization; for trial division, relation
deduplication, large-prime matching, matrix construction and the QS square root we found
no GPU prior art. Special-purpose sieving hardware~\cite{\geiselmanncite} is the
classical alternative.

\textbf{Concurrent work.} After our RSA-155 run, GPU adaptations of
CADO-NFS~\cite{cadonfs} factored RSA-260 (Lu, with the Devin agent: GPU lattice sieve,
Block Wiedemann and square root, CPU filtering)~\cite{lu2026rsa260} and the 270-digit
RSA-896 (Weis, with Claude)~\cite{weis2026rsa896}; both are GNFS.

\subsection{Quadratic-sieve records}
\label{sec:qsrecords}

The QS record is frequently misquoted as RSA-129 (1994; ppmpqs, ${\sim}600$
volunteers)~\cite{atkins1995squeamish}. That is the largest \emph{distributed
community} QS effort, and it ceased to be the size record in 2002; the record has since
moved three times, the last two by hobbyist effort: a 135-digit integer (2002; Leyland, Lenstra, Dodson, Muffett and Wagstaff,
MPQS with \emph{three} large primes~\cite{leyland2002threelp}); RSA-140 (June 2020, P.~Konsor, QS, ${\sim}6{,}000$
core-hours~\cite{konsor2020rsa140qs}, since improved to 1{,}791 by
Buhrow~\cite{buhrow2025rsa150qs}); RSA-150 (June 2025, B.~Buhrow, YAFU with the
\emph{three}-large-prime variant, 11{,}664 core-hours on AMD
EPYC~9174~\cite{buhrow2025rsa150qs,yafu}); and RSA-155 (July 2026, this work, SIQS
with a \emph{single} large prime, 64$\times$H100, 700.6~GPU-h).

Two features of this chain matter. Every QS size record since 2002 was set with a
\emph{multiple}-large-prime variation, whereas we use a single large prime
(\S\ref{sec:lp}): the record is set with the weakest large-prime variation of any entry
above it, which leaves the
multiple-large-prime variations as unclaimed headroom (\S\ref{sec:conclusion}). And our
margin over the prior QS record is 5 digits and 13 months.
General factoring records belong to the number field sieve: RSA-768~\cite{kleinjung2010rsa768},
RSA-250~\cite{boudot2020rsa250,factorrecords} and, concurrently with this work and on
GPUs, RSA-260 and RSA-896~\cite{lu2026rsa260,weis2026rsa896}. We did not factor a previously
unfactored number, and this is not a general factoring record. Buhrow's own assessment
corroborates the standing of the record we report against (``This last I think is the
largest factorization ever done by QS''~\cite{buhrow2025rsa150qs}); the comparison to
his RSA-150 is not like-for-like, so we run a head-to-head on the identical modulus
(\S\ref{sec:rsa150}).

\subsection{CPU implementations}

YAFU~\cite{yafu}, AVX-512 throughout, is the reference point for QS performance.
CADO-NFS~\cite{cadonfs} and msieve~\cite{msieve} are the standard NFS implementations;
CADO-NFS's development branch has since gained a CPU SIQS by Bouvier, a \texttt{siqs}
binary added December 2025 with Gray-code polynomial enumeration, two-level bucket
sieving and AVX2 sieve kernels --- the natural CPU-SIQS counterpart to the pipeline
described here, and in \S\ref{sec:cpubaseline} it runs ${\approx}2.4\times$ slower at
RSA-100 than YAFU's tuned single-large-prime SIQS. Bucket sieving,
deferring large-prime updates into per-subinterval buckets, is described by Aoki
and Ueda~\cite{aoki2004bucket} and Franke and Kleinjung~\cite{franke2005}; we adopt it
and implement it on GPU (\S\ref{sec:sieve}).
 \section{Design and Implementation}
\label{sec:design}

\subsection{Pipeline overview}

One orchestrator runs six stages: tuning, optional
autotuning (\S\ref{sec:autotune}), the sieve, matrix construction, linear algebra,
square root --- at v1.0.5, 70{,}458 lines and 145 unique \texttt{\_\_global\_\_}
kernels including the separately versioned \GFtwo{} Block Wiedemann submodule; C++20
host / CUDA~17 device, compute capabilities 7.5 through 12.0 (Turing to Blackwell,
Jetson Orin's 8.7 included).

\subsection{GPU residency, and its boundary}
\label{sec:residency}

Every stage executes in device kernels --- polynomial generation and its Gray-code
root update, both sieve kernels, trial division on a 1024-bit accumulator, relation
append and deduplication, device-side termination, large-prime matching, matrix construction,
all three Block Wiedemann stages with the Krylov sequence $S$ device-resident
and lingen's per-step loop \emph{including pivot selection}, and the square root. Five host residuals qualify ``fully on GPU''.

\paragraph{Residuals 1 and 2: setup and two materializations} Factor-base generation,
the Tonelli--Shanks square roots and the parameter heuristics are host code --- one-time
setup (13.5\,s measured at $F=5\times10^8$) against RSA-155's 10.78\,h sieve (${<}0.1\%$), and
deterministic, which cluster mode relies on (\S\ref{sec:cluster}). The relation set also
passes from sieve to matrix through a host structure-of-arrays container, and the matrix
from construction to Block Wiedemann through a host jagged-row container.

\paragraph{Residual 3: the hypercube walk} The choice of \emph{which} factor-base primes
enter $a$ is a host-side Gray-code sliding-window walk, run once per batch on the
critical path, producing $\textit{batch\_size}\times k$ 32-bit indices
(${\sim}2.4$\,KB) staged through pinned memory and copied asynchronously; all polynomial
arithmetic, the 512-bit $a$ and the $B_i$ basis included, is a device kernel, and under
graph capture the indices are pre-staged.

\paragraph{Residual 4: lingen initialization} Through v1.0.7 the Krylov sequence $S$
stays on device but is also copied to the host once per run, as $L$ separately
synchronized per-coefficient transfers ($O(L\,m\,n/64)$ words). On that copy one CPU
thread finds the initialization basis by incremental \GFtwo{} Gaussian elimination over
the columns of $S_0, S_1, \dots$, stopping when the rank reaches $m$ --- which fixes
$t_0$, so its cost scales with $t_0$, not $L$ --- and builds the initial generator and
degree vector, which are uploaded. v1.0.8 runs this step as a single-warp device kernel,
bit-identical to the host routine, and copies $S$ only on request; the RSA-150 and
RSA-155 solves used the host path. Every later step, pivot elimination included, runs
on device against the device-resident $S$.

\paragraph{Residual 5: multi-node large-prime matching} In cluster mode, large-prime
matching is deliberately relocated to a single coordinator CPU thread
(\S\ref{sec:cluster}, \S\ref{sec:rsa155}). Solo mode is GPU-resident; cluster
mode is not.

\subsection{The sieve}
\label{sec:sieve}

Sieving into shared-memory-sized subintervals caches well for small primes but
degenerates for large ones, whose roots are often loaded only to find the subinterval
holds no offset; \emph{bucket sieving}~\cite{aoki2004bucket,franke2005} defers those
writes by appending a 64-bit encoded entry (prime index, $\log_2 p$, offset) to the
bucket of the subinterval each update lands in. We call this step the
\emph{meta-sieve}: it sieves at the level of buckets rather than of offsets.
Sieve intervals being short, it amortizes each prime load by cycling through
\emph{several polynomials} of the current hypercube before advancing, a Gray-code
switch being a register-only add or subtract (\S\ref{sec:bg}); msieve's
SIQS~\cite{msieve} likewise nests a block of polynomials inside its loop over the
large primes.
That needs one bucket per (polynomial, subinterval) pair, and two
orthogonal subdivisions (of the subcube across thread blocks, of the subintervals
into sequential cycles) keep few enough active at once ($8\times16=128$ per thread
block) that their atomic fill counters live in \emph{shared} memory.

Per polynomial family the sieve is three kernels: \emph{initialization}, precomputing
$a^{-1}R_p$ and $a^{-1}B_i \bmod p$; the meta-sieve (scatter), filling buckets;
and \emph{sieve-and-scan} (gather), which per (polynomial, subinterval) performs the
queued writes, sieves the remaining small and intermediate primes directly, scans the
accumulator, and re-sieves in ``read'' mode to recover each candidate's factor list.
Gather parallelizes \emph{small} primes \emph{by offset} (unique write locations, no
races), \emph{intermediate} primes \emph{by prime} (one per thread, with atomics), and
\emph{large} primes entirely by the bucket; the small/intermediate boundary is a
nontrivial tuning parameter. The sieve claims the full opt-in dynamic shared-memory
budget, ${\sim}228$\,KB on an H100 rather than the 48\,KB static limit.

\subsection{The log accumulator: wrap, and a provable saturation gate}
\label{sec:accum}

The sieve accumulator is a byte array in shared memory and CUDA atomics are at least 32
bits wide, so updates go through a packed-byte \texttt{atomicAdd} on the containing
word, and the legacy accumulator wraps into the adjacent byte lane above 255.
At typical SIQS parameters that never fires; at RSA-155 it fires on exactly the values
that matter, a \emph{full} smooth accumulating $\log_2|Q/a|\approx277$, so the wrapped
accumulator lands \emph{below} threshold and the candidate is silently rejected --- an
estimated 98--99.5\% of full smooths lost, leaving an all-partial relation set and the
square-root failure of \S\ref{sec:limits}. RSA-140 accumulates ${\approx}249.9<255$ and
factored uneventfully; RSA-155 did not.

The fix is a dual-path accumulator selected at configuration time by an exact
three-point evaluation of the polynomial magnitude at $x\in\{-M,0,M-1\}$ mirroring the
device arithmetic. A \texttt{uint16} accumulator is always correct but \emph{halves}
the sieving block (131{,}072 $\to$ 65{,}536 entries on an H100) at a cost in
throughput, so we add a \emph{saturating} uint8 accumulator (\texttt{u8sat}) via a
read-clamp-CAS primitive, selected iff the config-time gate
\begin{equation}
  \textit{max\_target} \;=\; \textit{APV}_{\max} - \textit{threshold} \;\le\; 252
  \label{eq:u8satgate}
\end{equation}
holds. Under \eqref{eq:u8satgate}, $\textit{stored}=\min(255,\textit{sum})$ implies
$(\textit{stored} > \textit{target}) \Leftrightarrow (\textit{sum} > \textit{target})$
--- candidate selection is provably bit-identical to the uint16 path at one byte per
entry. RSA-155 runs at $\textit{max\_target}=249$ and RSA-150 at 240 (239 on the A100), so both
select \texttt{u8sat}; forcing it still honours the gate, never the wrapping path.
Measured, \texttt{u8sat}/uint16 is $1.56\times$ (H100, RSA-155) and $\ge1.65\times$
(A100, RSA-150).

\subsection{The zero-synchronization steady state}
\label{sec:zerosync}

We implemented {\em device-only} sieving in a strict sense: the solo batch sieve loop
contains no blocking host synchronization. Sieve and post-processing overlap
through a double-buffered, event-synchronized candidate handoff, the large-prime
pipeline runs on its own stream, and the batch boundary is an event record and wait.
Two details make the claim more than bookkeeping. First, termination is a
device decision: a yield-prediction kernel writes a mapped-pinned structure the host
polls \emph{without} synchronizing, so the loop decides it has enough relations
without draining the device. Second, large-prime telemetry is read \emph{deliberately
stale} from pinned memory, costing a measured ${\sim}1.2\%$ relation overshoot in
exchange for never stalling. On top of this, $N$ sieve batches can be captured into a
CUDA graph and replayed, capture being invalidated and re-taken adaptively when yield
drops.

Nsight Systems timelines of complete RSA-100 factorizations on an H100 (v1.0.5; the
v1.0.8 sieve takes 19.73\,s) measure the claim directly: across the sieve window the GPU is busy 99.9\% of the time with \emph{zero}
device-wide synchronizations, \emph{zero} event synchronizations and \emph{zero}
synchronous copies, and with graph capture off exactly \emph{one} blocking stream wait
survives, the final drain. Capture is not free everywhere: on that H100 kernels that overlap between batches
serialize between replays instead, each adding a stream wait, for $+3.3\%$ on total
run time --- so the discrete-GPU v1.0.8 benchmarks of \S\ref{sec:corescaling} and the
profiles of \S\ref{sec:roofline} run with capture off.

\subsection{On-device large-prime matching}
\label{sec:lp}

Large-prime matching is classically a host-side hash join; we do it on device, inside
the sieve loop, asynchronously. The structure is a slab hash table: a
directory of $2^B$ 64-bit words over payload rows of 16 eight-byte entries, locked per
bucket by \texttt{atomicCAS}, a full row rejecting gracefully instead of probing on.
Matching runs in three device stages --- directory fetch, vectorized 128-bit probe, and
commit, where matched pairs are re-checked against the actual large-prime values (tags
being truncated), merged and appended.

The acceptance gate tests only that the unfactored cofactor is below $L$. No primality
test is needed: every configuration here satisfies $L < f_{\max}^2$ (at RSA-100,
$10^{12} < (5.5\times10^6)^2$), and a cofactor below $f_{\max}^2$ surviving trial
division by every factor-base prime is necessarily prime.

\subsection{Matrix construction and \GFtwo{} linear algebra}
\label{sec:linalg}

Matrix construction is a device SoA$\to$CSR conversion with no reduction: a row-weight
kernel, a Thrust exclusive scan, a chunked fill kernel. An opt-in \emph{packed} GPU
preprocessor (singleton removal, weight-2 and Markowitz
merges, Montgomery-carried~\cite{montgomery1985} $\textit{sqrt\_Q}$ products) is a
negative result: it recovers ${\sim}34.7\%$ nontrivial dependencies against
legacy's ${\sim}56\%$ and collapses \GFtwo{} column diversity on high-large-prime data
(${\sim}42$K columns against a CPU preprocessor's ${\sim}143$K on the same input), so
RSA-155 was factored with \texttt{--matrix\_mode legacy}; its lasting value is
the genus-pinning diagnostic behind \S\ref{sec:limits}, and cycle filtering over $\ge3$
partials as in CADO-NFS~\cite{cadonfs} is the unimplemented path to a competitive one.

The Block Wiedemann engine is where the GPU pays. No vendor library provides \GFtwo{}
sparse linear algebra (cuSPARSE and cuBLAS are floating-point and integer only),
so the kernels are built from scratch. The engine
implements ten sparse and dense \GFtwo{} formats (M4RM, dense bit-slice,
warp-CSR, tiled COO and an unrolled variant, PForDelta and a bit-exact variant,
Delta-16, Golomb, ELLPACK) nine enabled by default. Auto-tuning is itself
GPU-resident, the matrix uploaded once and transposition, conversion and benchmarking
all on device, and a dynamic-programming solver assigns per-segment recipes,
different formats for different density regions of one matrix: at RSA-155, M4RM for the
56 densest rows, PForDelta for 4.81\,M, warp-CSR for 11.89\,M. Block widths default
to $m=n=256$, a default \S\ref{sec:rsa150} measures.

\paragraph{Stage 2: Block Berlekamp--Massey, on device, and quadratic}
The generating-polynomial stage computes $\Pi(x) \in
\mathbb{F}_2[x]^{(m+n)\times(m+n)}$ with $S(x)\cdot\Pi(x) \equiv 0 \bmod x^L$ by
Block Berlekamp--Massey on dense $(m+n)\times(m+n)$ \GFtwo{} matrices.
We do not implement Thom\'e's subquadratic divide-and-conquer
algorithm~\cite{thome2002lingen}: there is no recursive or middle-product layer in the
solver (it uses Karatsuba for polynomial multiplication, but no subquadratic lingen recursion). The
per-step loop (discrepancy, elimination, shift selection and the $F$ update, pivot
search included) runs entirely on device with no per-step host synchronization,
excepting, through v1.0.7, the host-side initialization at the stage boundary (\S\ref{sec:residency}). The
quadratic basecase is the engine's largest remaining cost: at RSA-155, Stage~2 took
4\,h\,14\,m\,18\,s, ${\approx}40\%$ of the solver (\S\ref{sec:rsa155}), so
implementing~\cite{thome2002lingen} is quantified headroom
(\S\ref{sec:conclusion}); the linear-algebra leg is 1.6\% of RSA-155's cost
despite it.

\subsection{Square root}

The square root is fully parallel over kernel vectors: a three-phase $X$ computation, a
parallel $Y$ computation (exponent accumulation, halving, batched exponentiation), a
batched GCD kernel, and a refinement kernel extracting complete factorizations by
coprime splitting. The CPU fallback triggers only when the GPU path finds no nontrivial
factor across \emph{all} solutions --- the failure mode of \S\ref{sec:limits}, where it
cannot help.

\subsection{The cluster layer}
\label{sec:cluster}

The multi-node design has one principle: the orchestrator's sieve loops are the
sole sieve implementation. Cluster mode injects a two-method \texttt{DataTap}
interface, \texttt{onBatchComplete()} and \texttt{shouldStop()}, at exactly two points
in the existing loop, under a contract that \texttt{onBatchComplete()} be non-blocking
($<50\,\mu$s). Three implementations serve all roles: \texttt{nullptr} in solo mode
(zero overhead), a direct channel for the coordinator's local sieve thread, and an async
network tap on workers, which memcpys into an SPSC ring and returns, leaving
serialization, heartbeats and stop-polling to a dedicated I/O thread. Tap overhead is
0.38\% at RSA-100 scale, and graph capture is unaffected because extraction
runs \emph{between} replays. Work is distributed as contiguous $a$-ranges weighted by SM
count $\times$ clock.

\paragraph{Workers run no large-prime matching} Remote workers and the coordinator's
local sieve thread all run with GPU large-prime matching disabled; matching happens on
the coordinator's CPU. That keeps the GPU sieve unblocked and gives the partials a single global meeting
point, a partial's matching partner being producible on any node; its cost is measured
in \S\ref{sec:rsa155}.
 \section{Evaluation}
\label{sec:eval}

\subsection{Methodology}

All times are total wall-clock, parameter tuning and, where enabled,
autotuning included; we report no aggregate excluding them, since that exclusion is not
consistently defined across runs. Named per-stage breakdowns reconcile to their totals
within 0.05\,s. Single-GPU runs use the platforms of Table~\ref{tab:rsa100}; cluster
runs use \hpcsite{} \texttt{gpu\_h100} partition, four H100-SXM5-94GB per node (132 SMs,
94\,GB HBM2e, 2396\,GB/s; SKU operator-confirmed), one GPU per MPI rank. Energy is not NVML \texttt{total\_energy\_consumption}, unsupported on this
SKU: every figure trapezoidally integrates power samples over actual timestamps, at the
cadence stated.

\subsection{Cross-hardware RSA-100}
\label{sec:rsa100}

\begin{table}[t]
\centering\footnotesize
\caption{RSA-100 (330-bit): total wall-clock, sieve and linear algebra (LA), same
source recompiled per architecture with per-architecture tuned parameters. ``Cores'' is
the CUDA-core count; the Jetson Orin Nano is a 25\,W integrated SoC.}
\label{tab:rsa100}
\setlength{\tabcolsep}{3pt}
\resizebox{\linewidth}{!}{%
\begin{tabular}{@{}llrrrrrr@{}}
\toprule
 & & & \textbf{v1.0.5} & \multicolumn{4}{c}{\textbf{v1.0.8}} \\
\cmidrule(lr){4-4}\cmidrule(lr){5-8}
\textbf{GPU} & \textbf{Arch.\ (CC)} & \textbf{Cores} & \textbf{Total (s)} &
\textbf{Total (s)} & \textbf{Sieve (s)} & \textbf{LA (s)} & \textbf{Sieve} \\
\midrule
H100 SXM5 (94\,GB) & Hopper (9.0)     & 16{,}896 & 51.12      & 29.18$^{\ddagger}$         & 19.73      & 8.97   & 67.6\% \\
RTX 5070 Ti        & Blackwell (12.0) & 8{,}960  & 85.19      & 51.01               & 36.97      & 13.65  & 72.5\% \\
A100 SXM4-40GB     & Ampere (8.0)     & 6{,}912  & 139.88     & 57.64               & 43.96      & 13.16  & 76.3\% \\
TITAN RTX          & Turing (7.5)     & 4{,}608  & 219.32     & 111.74              & 86.64      & 24.14  & 77.5\% \\
Jetson Orin Nano   & Ampere (8.7)     & 1{,}024  & 2{,}395.69 & 1{,}028.04$^{\ast}$ & 813.43     & 212.04 & 79.1\% \\
\bottomrule
\end{tabular}%
}

\smallskip
\begin{minipage}{\linewidth}\footnotesize
Each v1.0.8 row is the run with the median total of $n{=}3$ re-runs ($n{=}5$ on the
RTX~5070~Ti), each device at an operating point tuned for it, so a v1.0.8/v1.0.5 ratio is
a release-to-release comparison, not a binary-only one. Binary-only, at an
identical argument vector, v1.0.8 cuts total wall-clock by 18.9\% over a same-day
v1.0.7 control on the TITAN~RTX (137.85\,s) and by 20.8\% on the RTX~5070~Ti in one interleaved session
(64.23\,s to 50.87\,s on v1.0.8a, which the release matches within noise); the A100 row
runs a different $F$ and $M$ from its v1.0.7 tuple and gives no version ratio.
v1.0.8 was measured at RSA-100 only.
$^{\ddagger}$Same $F$ and $M$ as v1.0.7, kernel tuple re-searched for v1.0.8: 7.6\%
faster than the v1.0.7 tuple alternated in the same job, which itself runs 17.6\% faster
than on v1.0.7 (31.57 vs 38.30\,s; different job and node). $^{\ast}$Mean of $n{=}2$, on an older flag set than the other rows; at
the same argument vector and node, 35.9\% faster than v1.0.7 (1{,}602.64\,s). Remaining stages (H100 / RTX~5070~Ti / A100 / TITAN~RTX
/ Jetson): tuning 0.16 / 0.17 / 0.22 / 0.50 / 1.07\,s, matrix construction 0.02 / 0.06 / 0.02 / 0.09 /
0.13\,s, square root 0.29 / 0.16 / 0.29 / 0.36 / 1.38\,s.
All logged entries and the H100, A100 and TITAN~RTX
v1.0.5 totals trace to run logs retained in the artifact; the RTX~5070~Ti v1.0.5 total
is the median of $n{=}3$ logged re-runs of the v1.0.5 release (84.72\,s in our run
register), and the Jetson's comes from the run register, reproduced to 0.05\% by a
logged re-run (2{,}394.59\,s) after a first re-run's two kernel vectors gave only
trivial factors. The v1.0.5 column holds v1.0.4d--v1.0.5 builds (RSA-100 timings are
identical on v1.0.4 and v1.0.5, and the release leaves the 1.0.5a/b RSA-100 code path
unchanged; superseded earlier runs recorded 158.75\,s and 226.01\,s); it is the
generation the cost normalization of \S\ref{sec:benchfact} uses.
\end{minipage}
\end{table}

Two observations from Table~\ref{tab:rsa100} carry the rest of the paper. First,
on v1.0.8 the sieve is 68--79\% of total wall-clock on every platform and linear
algebra 21--31\% (v1.0.5: 77--91\% and 8--21\%): v1.0.7 and v1.0.8 target the sieve,
whose time fell 45--63\% on the discrete GPUs against 6--32\% for linear algebra; at
RSA-155 (v1.0.5) the split is 98.4\%/1.6\%.
Second, configuration is not held constant: $F$ ranges over 5.5--7\,M, $M$
over 131\,K--524\,K, batch size is 8 or 32, and in the v1.0.8 column every device
additionally runs a tuple of its own, so this is a portability measurement of
one source recompiled per architecture. To our
knowledge 51.0\,s on a consumer RTX~5070~Ti and 29.2\,s on
one H100 are the fastest end-to-end quadratic-sieve wall-clock times reported on any
platform.

\subsection{CPU baseline at RSA-100}
\label{sec:cpubaseline}

The CPU reference is our own controlled measurement, at RSA-100 on the same modulus,
of YAFU v3.1.9~\cite{yafu} (GCC 14, \texttt{znver5}, AVX-512/IFMA, msieve~1.38 block
Lanczos) and of CADO-NFS~\cite{cadonfs} (development \texttt{HEAD} of August 2026; its
SIQS and its GNFS), each on 96 threads of an AMD EPYC~9655 (Zen~5, no SMT),
\texttt{numactl}-pinned to one socket with the affinity verified live, on exclusive
nodes of \hpcsite{} CPU cluster. Both sides are device-only, and both CPU wall-clocks
bracket the whole run as ours do. Package energy integrates 5\,s AMD HSMP per-socket
\texttt{hwmon} samples, independently reproduced; that sensor is uncalibrated vendor
telemetry with no second instrument, so row-to-row ratios are solid while \emph{absolute}
watt-hours are only indicative. On the GPU side, the H100 row is the median-total repetition of our v1.0.8
benchmark (job~3474713), time and energy from the same run: 200\,ms \texttt{nvidia-smi}
board power integrated over actual timestamps gives 3.38\,Wh per factorization
(3.37--3.41 over $n{=}3$), mean board power 406\,W against a 70.04\,W idle baseline. The
H100's v1.0.7 tuple, alternated in the same job, is 7.6\% slower at 1.7\% less energy.
Host brackets have one-second resolution, bounding each integral within $-1.5/{+}1.9\%$,
and the job shared its node, which
leaves the allocated card's own sensor unaffected. The on-card counter being unsupported
on this SKU, accuracy is dominated by the NVML power specification (${\sim}\pm5\%$). The
same instrument on v1.0.5, in a dedicated exclusive $n{=}4$ re-run reproducing the
51.12\,s best time within 0.21\%, gave 4.876\,Wh at 340.8\,W.

\begin{table}[t]
\centering\footnotesize
\caption{RSA-100 on one H100 against our measured CPU runs, same modulus throughout.
\textbf{LP} is large primes admitted per relation; wall is total wall-clock. Energy is a
power-sampler integral on both sides: CPU \textbf{package} (both sockets; loaded socket
alone in parentheses) from 5\,s samples, GPU \textbf{board} from 200\,ms samples;
``$\div$H100'' divides the row's package energy by the H100's 3.38\,Wh (v1.0.8).}
\label{tab:cpubase}
\setlength{\tabcolsep}{3pt}
\resizebox{\linewidth}{!}{%
\begin{tabular}{@{}llcrrr@{}}
\toprule
\textbf{Software} & \textbf{Hardware} & \textbf{LP} & \textbf{Wall (s)} &
\textbf{Energy (Wh)} & \textbf{$\div$H100} \\
\midrule
YAFU, default FB$^{a}$      & EPYC 9655, 96\,c.$^{d}$ & 1 & 157.67 & 20.59 (16.56) & 6.1$\times$ \\
\textbf{YAFU, tuned FB}$^{a}$ & EPYC 9655, 96\,c.$^{d}$ & 1 & \textbf{122.10} & \textbf{15.66} (12.50) & \textbf{4.6}$\times$ \\
\textbf{YAFU, default FB}$^{b}$ & EPYC 9655, 96\,c.$^{d}$ & 2 & \textbf{106.30} & \textbf{12.67} (10.14) & \textbf{3.7}$\times$ \\
YAFU, default FB$^{b}$      & EPYC 9655, 96\,c.$^{d}$ & 3 & \multicolumn{3}{l}{\emph{did not factor} (see text)} \\
CADO-NFS SIQS, tuned$^{c}$  & EPYC 9655, 96\,c. & $\le$2 & 291.85 & 40.15 (31.53) & 11.9$\times$ \\
CADO-NFS GNFS, default$^{c}$ & EPYC 9655, 96\,c. & --- & 262.76 & 22.79 (15.07) & 6.7$\times$ \\
\midrule
\textbf{This work} (v1.0.8) & 1$\times$ H100 SXM5 & 1 & \textbf{29.18} & \textbf{3.38}$^{\dagger}$ & 1 \\
\bottomrule
\end{tabular}%
}

\smallskip
\begin{minipage}{\linewidth}\footnotesize
$^{a}$\,One job, one node; the default-FB row is that job's pooled $n{=}3$ control, so
the two 1-LP rows differ only in factor base; it recorded no binary hash, and the
campaign's safer factor base reproduces the tuned optimum to 0.08\%.
$^{b}$\,Second job and node.
$^{c}$\,Third job and node, $n{=}3$ per row (spread 0.9\% and 1.5\%); wall is the run
bracket, CADO's own elapsed is 0.2--0.3\% lower. Upstream ships no QS parameter file
above 90 digits, so we wrote one for 100 and tuned it (factor-base bound, large-prime
bound, relation target; worth $-4.4\%$ wall at 96 threads); it admits cofactors up to
$2^{48}$, i.e.\ up to two large primes. GNFS is CADO's shipped \texttt{c100} default.
$^{d}$\,96 threads is this socket's measured optimum: at the
same configuration 48 and 24 threads cost $+56\%$ and $+193\%$ wall-clock and
$1.6\times$ and $2.5\times$ package energy (job~2965631, socket-0 pinned, exclusive,
same binary path and flags as the rows above; $96/48/24$\,T $=
122.10/190.92/358.10$\,s tuned 1-LP, 96\,T being row two). Beyond one socket untested.
$^{\dagger}$\,v1.0.8, median-total run of $n{=}3$ (job~3474713). Against the loaded socket alone
the ratios are 4.9$\times$, 3.7$\times$, 3.0$\times$, 9.3$\times$ and 4.5$\times$ (top
to bottom). On v1.0.5 (51.12\,s, 4.88\,Wh) every energy ratio is $1.44\times$ and every
wall-clock ratio $1.75\times$ smaller.
\end{minipage}
\end{table}

YAFU has no 1-LP mode: we force one with \texttt{-siqsMFBD 1.0}, which empties its
double-large-prime window. Its default factor base is tuned for two, so we swept
\texttt{-siqsB} over 20 legs and quote the optimum (208{,}400 primes), worth
$-22.6\%$ wall-clock and $-24.0\%$ energy against the same job's control. The 1-LP rows
are the algorithmic match, their cycle graph having no cycle longer than two, exactly
the relation algebra our pipeline produces; the 2-LP row is YAFU's best showing,
faster and lower-energy.
The 3-LP run did not factor RSA-100: a defect in YAFU's own
\texttt{check\_relation} rejected all 402{,}235 relations. The best published CPU
quadratic-sieve cost at RSA-100 --- 0.96 core-hours with YAFU's three-large-prime
variant on an EPYC~9174, reported in core-hours only~\cite{buhrow2025rsa150qs} --- is thus a 3-LP result we could not reproduce; our measured best, 2.835 core-hours, is ${\approx}3\times$ higher: part large-prime
variant, part tuning, part clock and per-core L3, QS sieving being cache-resident.

Three conclusions follow, in decreasing order of what the data support.
First, YAFU remains the fastest CPU quadratic sieve and its rows are therefore
the CPU-QS bar: its tuned 1-LP run is $2.39\times$ faster than CADO's tuned SIQS
(122.10 vs 291.85\,s) at $2.56\times$ less package energy. Second, one H100 on v1.0.8 is
$4.2\times$ faster than the fastest CPU quadratic sieve, $10.0\times$ faster than
CADO's SIQS and $9.0\times$ faster than its GNFS, at the package-energy ratios of
Table~\ref{tab:cpubase}'s last column; against YAFU's best 2-LP setting the margins are
$3.6\times$ in wall-clock and $3.7\times$ in energy (v1.0.5: $2.4$, $5.7$, $5.1$, $2.1$
and $2.6\times$). Every energy ratio compares
\emph{GPU board} with \emph{CPU package}; neither is node-level, and the instruments
differ. Third, RSA-100 sits at the upper end of the range in
which the quadratic sieve is competitive with the number field sieve: CADO's shipped
GNFS already beats its own tuned SIQS here (262.76 vs 291.85\,s, at $1.76\times$ less
energy), and is not even at its optimum, running \emph{faster} at 174.5\,s on 16
socket-confined threads of a like node. That is not a crossover measurement (one
implementation, one modulus, a CPU SIQS still under development) but it is same-machine
evidence for the ${\sim}$100--110-digit takeover, above which a quadratic-sieve run is
of systems interest and is no longer the cheapest way to factor.

\subsection{Scaling in CUDA cores}
\label{sec:corescaling}

\begin{figure}[t]
\centering
\begin{tikzpicture}
\begin{axis}[
  width=0.95\linewidth, height=3.8cm,
  xlabel={CUDA cores}, ylabel={RSA-100 throughput $1/t$ \; [$10^{-3}\,\mathrm{s}^{-1}$]},
  xmin=0, xmax=18200, ymin=-1.5, ymax=38,
  xtick={0,4000,8000,12000,16000},
  scaled x ticks=false,
  xticklabels={0,{4{,}000},{8{,}000},{12{,}000},{16{,}000}},
  grid=both, major grid style={line width=0.2pt, draw=gray!30},
  tick label style={font=\footnotesize}, label style={font=\small},
]
\addplot[thick, black!45, dashed, domain=0:18200, samples=2] {0.0012189*x -0.6369};
\addplot[only marks, mark=o, mark size=2.6pt, black!55] coordinates {
  (4608, 4.560) (6912, 7.149) (8960, 11.738) (16896, 19.562)
};
\addplot[thick, blue, domain=0:18200, samples=2] {0.0019492*x +1.8300};
\addplot[only marks, mark=*, mark size=2.6pt, red] coordinates {
  (4608, 8.949) (6912, 17.349) (8960, 19.604) (16896, 34.270)
};
\addplot[only marks, mark=square*, mark size=2.6pt, black!60] coordinates {
  (1024, 0.973)
};
\node[anchor=south east, font=\scriptsize] at (axis cs:4608,8.949) {TITAN RTX};
\node[anchor=south, font=\scriptsize] at (axis cs:6912,17.349) {A100};
\node[anchor=north west, font=\scriptsize] at (axis cs:8960,19.604) {RTX\,5070\,Ti};
\node[anchor=east, font=\scriptsize] at (axis cs:16896,34.270) {H100};
\node[anchor=west, font=\scriptsize] at (axis cs:1300,0.973) {Jetson};
\end{axis}
\end{tikzpicture}
\caption{RSA-100 throughput ($1/t$, total wall-clock) versus CUDA-core count, same
source recompiled per architecture. Hollow circles and the dashed line are the v1.0.5
anchors of Table~\ref{tab:rsa100} and their with-intercept fit ($R^2=0.978$); filled
circles and the solid line are the v1.0.8 benchmark (medians) and its fit
($R^2=0.976$).
The 25\,W Jetson Orin Nano (square, $n{=}2$; an integrated SoC on shared
LPDDR5) is plotted but excluded from both fits.}
\label{fig:corescaling}
\end{figure}
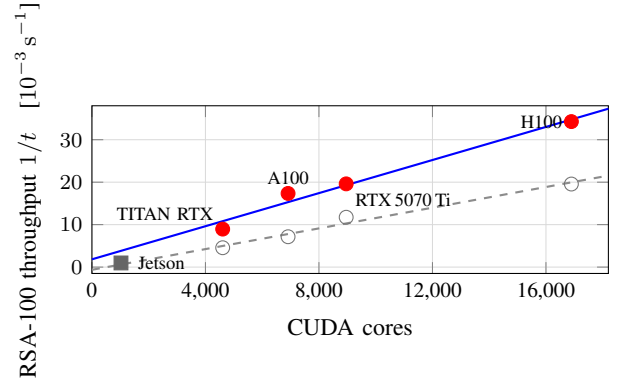

Across four discrete GPUs spanning four architectures and a $3.7\times$ range in CUDA
cores, RSA-100 throughput is close to proportional to core count on both releases
(Fig.~\ref{fig:corescaling}): the with-intercept linear fits are $y = 0.0012189\,x -
0.6369$ at $R^2 = 0.978$ on the v1.0.5 anchors and $y = 0.0019492\,x + 1.8300$ at
$R^2 = 0.976$ on the v1.0.8 benchmark, a $59.9\%$ steeper slope, with per-core
throughput constant to within $1.32\times$ and $1.29\times$ respectively. The intercept
is small and its sign is not stable across the two series ($-0.64$ against $+1.83$), so
at $n=4$ it carries no information; we plot both fits as computed and do not force
either through the origin. Memory bandwidth, which
varies ${\sim}3.6\times$ across the same parts (672\,GB/s to the 2396\,GB/s of this
H100 SKU --- the bus width and memory clock the device itself reports, against a nominal
2446\,GB/s in the operators' system description), remains a weaker predictor on
both ($R^2 = 0.64$ and $0.81$). \S\ref{sec:roofline} shows why: the sieve kernels
execute no floating-point work and draw about half of peak DRAM bandwidth or less (at
most 50.1\%) at every operating point profiled, even where v1.0.8 puts scatter below the ridge.

The residuals track operating point, not architecture: $-8\%$ (TITAN~RTX, A100) to
$+14\%$ (RTX~5070~Ti) on v1.0.5, and $-17\%$ / $+13\%$ / $+2\%$ / $-1\%$ (TITAN~RTX /
A100 / RTX~5070~Ti / H100) on v1.0.8, where the A100 (re-tuned, $M=524$\,K) and the
H100 run new tuples, the TITAN~RTX and RTX~5070~Ti their v1.0.7 ones.

Three remarks. The Jetson is
excluded from both fits: it is a 25\,W edge SoC with ${\sim}102$\,GB/s of shared LPDDR5
against 672--2396\,GB/s of discrete memory, and its point, 75\% below the v1.0.8 fit, would widen the
v1.0.8 per-core spread from $1.29\times$ to $2.64\times$, so the ``constant per-core
throughput'' claim fails with it (its point is also the weakest, $n{=}2$ on an older
flag set).
What the fit establishes is peak compute, not core count \emph{specifically}:
on v1.0.5 cores and cores$\times$clock are indistinguishable ($R^2 = 0.978$ for
both), clock varying only $1.4\times$ across these parts against $3.7\times$ in
cores, and while cores alone does separate on v1.0.8 ($0.976$ vs.\ $0.876$), reading a
mechanism out of that would over-interpret four points --- with 2 residual
degrees of freedom, even $R^2 = 0.976$ is only suggestive. Finally,
between the two series each card changes both binary \emph{and} tuned operating
point, so the gap between them is not a binary-only speedup.

\subsection{Benchmark factorizations and scaling in \texorpdfstring{$N$}{N}}
\label{sec:benchfact}

\begin{table}[t]
\centering\footnotesize
\caption{\swname{} benchmark factorizations. GPU-hours sum the
sieving GPUs and the single-GPU linear-algebra leg; \textbf{LP\,\%} is the share of
relations that are large-prime combines.}
\label{tab:benchfact}
\setlength{\tabcolsep}{4pt}
\begin{tabular}{@{}llrr@{}}
\toprule
\textbf{Modulus} & \textbf{Hardware} & \textbf{GPU-h} & \textbf{LP\,\%} \\
\midrule
RSA-100 & 1$\times$ H100             & 0.014  & 52.1 \\
RSA-110 & 1$\times$ RTX 5070 Ti      & 0.238  & 54.7 \\
RSA-120 & 1$\times$ RTX 5070 Ti      & 2.95   & 48.7 \\
RSA-130 & 8$\times$ A100 (2 nodes)   & 25.75$^{a}$ & 51.8 \\
RSA-140 & 16$\times$ H100 (4 nodes)  & 103.9$^{b}$ & 56.8 \\
RSA-150 & 64$\times$ H100 (16 nodes) & 302.9  & 44.9 \\
RSA-155 & 64$\times$ H100 (16 nodes) & 700.6  & 41.4 \\
\bottomrule
\end{tabular}

\smallskip
\begin{minipage}{\linewidth}\footnotesize
$^{a}$\,A100-GPU-hours; all others H100-GPU-hours.
$^{b}$\,Ran at deduplication safety factor 1.4 (40\% oversieve) versus 1.05 elsewhere; a
like-for-like target implies $\approx$78\,GPU-h. Energy instrumentation exists for
RSA-150 and RSA-155 only.
\end{minipage}
\end{table}

Table~\ref{tab:benchfact} lists them, RSA-100 to RSA-120 on a single device and
RSA-130 up on clusters. The v1.0.6--v1.0.8 sieve optimizations target the narrow 8-bit
accumulator path; RSA-150 and RSA-155 run the wide saturating-8-bit path, whose gather
and scatter kernels are byte-identical from v1.0.5 through v1.0.7 (v1.0.8 shares the
rewritten scatter and edits the wide gather, unmeasured there) and whose parameters
the new tuner declines to search. We therefore report those v1.0.5 results unchanged and make
no version-related projection; RSA-130 and RSA-140 are v1.0.3b and v1.0.4d
measurements, and we project nothing for them either.
Normalizing the non-H100 rows to H100-GPU-hours by the \emph{v1.0.5} RSA-100 ratios of Table~\ref{tab:rsa100}
(RTX~5070~Ti $\div1.663$, A100 $\div2.739$) and fitting $\ln(\text{cost}) =
e\ln\LN{1/2}{1} + \ln C$ over the seven factorizations gives $e = 1.04$ at $R^2 = 0.993$; the
four direct-H100 points alone (RSA-100, 140, 150, 155) give $e = 1.05$.
Since heuristic QS complexity is $\LN{1/2}{1}^{1+o(1)}$, a finite-range fit folds the
$o(1)$ term, per-run tuning and hardware saturation into its exponent, and the
residuals are dominated by tuning (RSA-120 and RSA-140 sit high, RSA-150 and RSA-155
22\% and 27\% below): the exponent summarizes these seven runs, and is not an
asymptotic statement about the algorithm. We do not compare it with exponents fitted
to CPU results --- the large-prime variants differ algorithmically, and a finite-range
exponent cannot separate variant, tuning and hardware.

\subsection{RSA-150: the same-modulus head-to-head}
\label{sec:rsa150}

The most informative comparison available to this work is RSA-150, the identical
modulus on which Buhrow's YAFU run set the prior QS record at 11{,}664 core-hours (about four
calendar days on 128 cores) with the three-large-prime variant on AMD EPYC~9174 processors~\cite{buhrow2025rsa150qs}
(its block-Lanczos step under an hour on 32 threads by his log, negligible against
the sieve as in our runs). It is a comparison and not a record: RSA-150's factors have been public since the
number field sieve produced them in 2004~\cite{factorrecords}, and the QS record we
claim is RSA-155 (\S\ref{sec:rsa155}). What it buys is an identical $N$ measured in
both directions.

We factored it for 302.86 H100-GPU-h and 107.5\,kWh: a 64$\times$H100 sieve of
295.18 GPU-h (4.612\,h wall) producing 17{,}269{,}643 relations (44.9\% large-prime
combines), then a single-H100 \GFtwo{} solve of 7.68 GPU-h (7.675\,h wall) on a
$16.70\text{M}\times15.66$M matrix with $6.9\times10^{8}$ nonzeros.
Sieving is 97.5\% of the GPU-hours and 98.3\% of the energy, as at RSA-155, the energy
again a power-sampler integral. The two legs were \emph{separate scheduler jobs},
48\,h of queue apart.
Against RSA-155 the \emph{sieve} cost ratio is $2.337\times$, inside the
$2.32$--$2.39$ band the asymptotics predict. The prediction we committed to \emph{before}
the run, from the then-six-point fit of \S\ref{sec:benchfact}, was $\approx435$ GPU-h,
band 256--616, pass criterion 300--450: the band held and the point estimate
over-predicted by 44\%, the fit absorbing an average per-run inefficiency a tuned run
does not pay. The sieve's own pre-run ETA, 292 GPU-h, came in at $+1.1\%$.

Against Buhrow the ratio is 38.5 core-hours per GPU-hour --- a
resource-replacement ratio, not a speedup, the two units not being commensurable. Three caveats all cut against us. He runs \emph{three} large primes to our \emph{one},
an algorithmic advantage that grows with $N$; he disclaims his own tuning (``probably
far from optimal''~\cite{buhrow2025rsa150qs}), so 11{,}664 core-hours is not a tuned
lower bound for YAFU, a log-log fit through his six data points sitting 12\%
\emph{below} what he measured; and the energy comparison does not survive the choice
of CPU part, a rated-TDP model running $2.2\times$ in our favour on the EPYC~9174F (his ``9174'')
and $0.6\times$ on a core-dense EPYC of the same generation, consequently we draw no
energy conclusion.

\paragraph{Block width dominates the linear algebra}
Three linear-algebra runs consumed the same saved relation set and built a
bit-identical matrix --- 690{,}586{,}752 nonzeros, same binary and solver, same flags
but $m$ and $n$ --- so block width is a controlled variable. At $m=n=256$ the production solve took 7.68 GPU-h; at $m=n=128$ it took 20.60
GPU-h, $2.69\times$ more, a loss of 12.92 GPU-h and 1.91\,kWh of which 59.7\%
is lingen. Halving the width doubles the Krylov sequence length
($130{,}596\to261{,}000$) and every degree that follows, and leaves each step \emph{less
efficient} besides, a $\dim=256$ step no longer filling an H100 (151\,W against 236\,W);
the projection that sent us to $n=128$ was refuted in reconstruction, whose per-step
cost scales with $n$ and step count as $1/n$. Time was not the worst of it: the narrow block found \emph{three} independent solutions where
the wide one found 56, keeping 24 at the cap, and at our measured 45.8\% nontrivial rate three
solutions carry a $(1-0.458)^3 = 16\%$ chance of yielding no factor at all after 20.6
GPU-h, against ${\approx}4\times10^{-7}$ at $m=n=256$. Width is bracketed from above, so
this is an optimum: a probe at $m=n=512$ ran lingen
$3.98\times$ slower ($\dim=1024$ falls off the specialized dense \GFtwo{} matmul paths)
and then exhausted memory in reconstruction after 18.24\,h, putting machine time
actually spent on RSA-150 at 323.45 GPU-h. One variable is uncontrolled: the
three runs took different non-exclusive nodes on different days, a drift bounded by the
block-independent SpMM autotune stage, which varied by 3.8\% across them --- two orders
of magnitude below the effect claimed. Lingen, quadratic Block Berlekamp--Massey, is
55.3\% of the winning solve (4\,h\,14\,m, as at RSA-155: its cost depends only on the
sequence length and block width, equal at both), the largest single piece of evidence for the headroom of
\S\ref{sec:conclusion}.

\subsection{QS RSA-155 factorization}
\label{sec:rsa155}

On 14 July 2026 the pipeline factored the 512-bit RSA-155 challenge modulus end to end
into two 78-digit primes, product-verified on device and re-verified host-side.

\paragraph{Cost} A 64$\times$H100 sieve across 16 nodes ($1{,}081{,}344$ CUDA cores)
ran 10.78\,h wall (689.74\,GPU-h, 239.78\,kWh), feeding a single-H100 \GFtwo{} solve
of 10.90\,h (10.90\,GPU-h, 2.48\,kWh): 700.64\,GPU-h and 242.26\,kWh,
${\approx}24$\,h from sieve start to factors across two jobs, the sieve 98.4\% of GPU-hours and 99.0\% of
energy.

\paragraph{Energy}
\label{sec:energy}
The RSA-155 energy figures are trapezoidal integrals over actual timestamps of 10\,s
power samples (247{,}857 across the 64 sieve ranks, 3{,}923 for the linear-algebra
leg). The sieve ran at a mean of 347.6\,W per H100 against a 700\,W cap
(per-GPU means 333--362\,W): the sieve is not power-bound.
Energy efficiency is not claimed as a contribution: the only CPU-versus-GPU energy
comparison is the measured one at RSA-100 (Table~\ref{tab:cpubase}), and we model none
at larger $N$.

\paragraph{Configuration} $F = 6\times10^8$ (15{,}666{,}200 factor-base primes), $M =
8{,}388{,}608$, $L = 8\times10^{13}$, the \texttt{u8sat} wide accumulator on all 64
ranks (gate \eqref{eq:u8satgate} holds at $\textit{max\_target}=249$), batch size 16
at CUDA-graph unroll 4, autotuning off (\S\ref{sec:autotune}), matrix stage
\texttt{legacy} (\S\ref{sec:linalg}).

\paragraph{Yield} 17{,}272{,}258 relations: 10{,}121{,}253 full smooths (58.6\%) and
7{,}152{,}282 large-prime combines (41.4\%), from 114{,}590{,}034 raw partials and
100{,}279{,}366 witnesses (6.7\% yield); aggregate 446.3 relations/s over 38{,}703\,s.
The 41.4\% large-prime share is structural to $F=6\times10^8$, not a tuning miss: the
run's telemetry confirms the birthday model ($\textit{combines} \propto
\textit{witnesses}^{1.988}$, $R^2=0.99999$), and the Mertens window at this $F$ is
such that raising $L$ from 80\,T to 200\,T would lift it only to ${\sim}44\%$.

\paragraph{The coordinator is not the limit} In cluster mode all large-prime matching
runs on one coordinator CPU thread (\S\ref{sec:cluster}). An estimate inferred from this
run put that thread's ceiling near the fleet's ${\approx}2{,}600$ witnesses/s, and a
64$\to$108-GPU follow-up was declined on it; direct measurement refutes it. At RSA-150
the same thread absorbed ${\approx}9{,}700$ witnesses/s while asleep 97.0\% of the run,
and the coordinator process peaked at 110\,GiB of its 512\,GiB allocation. 64 GPUs is a single
operating point, not a measured limit; we report no GPU-count sweep.

\paragraph{Linear algebra} The matrix was $16{,}700{,}000 \times 15{,}666{,}202$ with
683{,}983{,}246 nonzeros (41.0 per row), the row cap keeping the padded square
dimension below $2^{24}-1$ for tiled-COO admissibility at no cost. Block Wiedemann ran
at $m=n=256$ and sequence length 130{,}596, 10\,h\,26\,m in total (Krylov 1\,h\,23\,m,
lingen 4\,h\,14\,m, reconstruction 4\,h\,49\,m), matrix construction 1.57\,s and the
square root 15.22\,s. Of
24 independent kernel solutions, all 24 verified $X^2 \equiv Y^2$ and
11 (45.8\%) yielded the nontrivial split.

\subsection{Bottleneck analysis}
\label{sec:roofline}

\begin{table*}[t]
\centering
\scriptsize
\setlength{\tabcolsep}{1.2pt}
\caption{Bottleneck analysis of the two sieve kernels (Sec.~\ref{sec:sieve}) at RSA-100.
Columns~1--2 are one A100 on v1.0.6 binaries under power-of-two (v1.0.5-era) and
SM-aligned (v1.0.6) launch geometry; columns~3--4 the same device at the same SM-aligned pin
on the released v1.0.7 and v1.0.8; column~5 an H100 on the released v1.0.8, power-of-two.}
\label{tab:roofline}
\resizebox{\textwidth}{!}{%
\begin{tabular}{lrrrrr@{\hspace{5pt}}rrrrr}
\toprule
 & \multicolumn{5}{c}{SCATTER} & \multicolumn{5}{c}{GATHER} \\
\cmidrule(lr){2-6} \cmidrule(lr){7-11}
Metric & \shortstack{A100\\v1.0.5-geom} & \shortstack{A100\\v1.0.6-geom} & \shortstack{A100\\v1.0.7} & \shortstack{A100\\v1.0.8} & \shortstack{H100\\v1.0.8} & \shortstack{A100\\v1.0.5-geom} & \shortstack{A100\\v1.0.6-geom} & \shortstack{A100\\v1.0.7} & \shortstack{A100\\v1.0.8} & \shortstack{H100\\v1.0.8} \\
\midrule
GPU-time share (\%) & 34.1 & 26.0 & 39.2 & 47.0 & 44.0 & 48.7 & 55.6 & 52.8 & 41.4 & 43.8 \\
Duration/launch (ms) & 3.391 & 3.689 & 3.668 & 2.836 & 2.366 & 4.859 & 7.847 & 5.124 & 2.667 & 2.470 \\
Grid (blocks) & 64 & 108 & 108 & 108 & 256 & 512 & 864 & 864 & 864 & 256 \\
DRAM (GB/s) & 259 (16.6\%) & 403 (25.9\%) & 405 (26.1\%) & 526 (33.8\%) & 1180 (49.3\%) & 184 (11.8\%) & 191 (12.3\%) & 293 (18.9\%) & 572 (36.8\%) & 1159 (48.4\%) \\
Arith.\ int.\ (int-op/B) & 18.09 & 18.02 & 17.80 & 8.48 & 6.01 & 14.41 & 14.50 & 13.31 & 8.31 & 8.89 \\
INT32 (\% of peak) & 48.0 & 74.5 & 74.0 & 45.8 & 42.4 & 27.2 & 28.4 & 40.1 & 48.7 & 61.6 \\
L2 hit (\%) & 99.19 & 99.32 & 99.31 & 99.16 & 98.94 & 28.35 & 28.35 & 29.01 & 28.67 & 3.95 \\
Top stall (\% of budget) & \texttt{not\_sel} 25.5 & \texttt{not\_sel} 22.7 & \texttt{not\_sel} 23.0 & \texttt{long\_sb} 33.9 & \texttt{long\_sb} 43.3 & \texttt{barrier} 52.0 & \texttt{barrier} 52.2 & \texttt{barrier} 34.2 & \texttt{barrier} 32.8 & \texttt{long\_sb} 25.8 \\
\bottomrule
\multicolumn{11}{p{0.97\textwidth}}{\scriptsize ncu means over 10 steady-state launches per kernel (--launch-skip 20 --launch-count 10, --clock-control none, --cache-control none, --cuda\_graph\_unroll 0 so streams are serialized); share of GPU time from nsys (A100 v1.0.5-geom and v1.0.6-geom: full pipeline, --cuda\_graph\_unroll 4; A100 v1.0.7, v1.0.8 and H100 v1.0.8: 200-batch census, --cuda\_graph\_unroll 0 -- so that row is NOT comparable across columns). Percent-of-peak is against each column's own device's driver-reported peaks (A100-SXM4-40GB: 1555.2 GB/s DRAM, 9.75 T int-op/s; H100: 2395.9 GB/s DRAM, 16.73 T int-op/s). Arithmetic intensity is integer thread-ops (smsp\_\_sass\_thread\_inst\_executed\_op\_integer\_pred\_on.sum) per DRAM byte (dram\_\_bytes.sum) -- the instruction roofline of~\cite{ding2019roofline}; both kernels execute zero floating-point instructions on every metric collected. Warp-stall entries name the leading reason (long\_sb = long\_scoreboard, not\_sel = not\_selected) and its share of the stall budget (ncu's Warp Cycles Per Issued Instruction less the one 'selected' cycle, i.e. every stall reason whether or not it was exported). Achieved / theoretical occupancy is 50.0 / 50 in every column but v1.0.8 scatter, where two blocks fit per SM: 50.1 / 100 on the A100 (the 108-block grid fills half), 92.6 / 100 on the H100. A100 jobs 34116294, 34117902, 34507954 arm \texttt{std7} (its in-job v1.0.6 control reproduces column~2) and 34564985 arm \texttt{std8} (its in-job v1.0.7 control reproduces column~3's gather instruction count exactly); H100 job 3474714 arm \texttt{std8} (its in-job v1.0.7 control reproduces job 3170077). Columns~3--5 are each device's own v1.0.6 pin, \emph{not} the tuned tuples behind Table~\ref{tab:rsa100}.} \\
\end{tabular}%
}
\end{table*}

Why does throughput track peak compute (\S\ref{sec:corescaling})? We profiled RSA-100
with Nsight Compute and Nsight Systems (methodology in Table~\ref{tab:roofline}'s note;
a cold-cache control differs by $\le0.11\%$, the two tools agreeing to $\le0.4\%$ on v1.0.6): the
A100 under v1.0.6 (both launch geometries), v1.0.7 and v1.0.8, the H100 under v1.0.7
and v1.0.8 --- each version comparison made within one job at the same pin.

The two sieve kernels of \S\ref{sec:sieve} are 82.8\% of traced GPU kernel
time in the A100 full-pipeline trace and 88--92\% of a sieve-only census,
against ${\approx}12\%$ for linear algebra and 0.12\% for trial division.
Both are \emph{floating-point-free} on every metric collected, so
Table~\ref{tab:roofline} places them on the \emph{instruction}
roofline~\cite{ding2019roofline}: at the pins, 17.8 and 13.3 integer operations per DRAM
byte on the A100 under v1.0.7 and 8.5 and 8.3 under v1.0.8, against a machine balance of
6.3; on the H100 12.1 and 15.0 under v1.0.7 and 6.0 and 8.9 under v1.0.8, against 7.0.
Every column but H100 v1.0.8 scatter sits 1.3--2.9$\times$ compute-side of the integer
ridge; v1.0.8's smaller instruction stream moves scatter below it on the H100
(0.86$\times$ at the pin, 0.95$\times$ at its tuned tuple) and at the A100's tuned tuple
(0.95$\times$). Yet mean DRAM draw never exceeds 50.1\% of peak at any point profiled, and
\texttt{long\_scoreboard} leads each sub-ridge scatter (42--44\% of the stall
budget): the ridge is crossed by removing integer work, not by
saturating DRAM. Neither kernel is bandwidth-bound, though the handoff between
them is a full DRAM round-trip by design (scatter writes 0.9--1.5\,GB of buckets per
launch at 98--99\% L2 hit and gather reads them back at 28--29\% L2 on the A100, a working set
${\sim}35\times$ the L2): those scattered reads cost latency.

Through v1.0.7 scatter is integer-ALU-issue-bound, its instruction stream
identical across v1.0.6 and v1.0.7: ALU pipe at 74--79\% of peak, its leading warp states
describing warps ready but unpicked or throttled at that pipe, the 82\% shared-memory
atomic bank-conflict rate secondary. v1.0.8 cuts that stream ($-48\%$ and $-39\%$ thread
instructions at the A100 and H100 pins) and registers from 40 to 32, so two blocks fit per
SM; the kernel turns latency-bound (\texttt{long\_scoreboard} leads at 34\% and 43\%)
and its profiled duration drops only 23\% and 8\%, the bucket-write stream and shared atomics
being unchanged. Gather's limiter is neither version- nor
tuple-invariant. Under v1.0.6 it is barrier-bound on both devices ---
\texttt{\_\_syncthreads()} is 52\% of the A100 stall budget, 25\% of the H100's --- because
registers and 90.8\,KB of dynamic shared memory \emph{independently} cap it at one
1024-thread block per SM, nothing co-resident hiding the barrier.
v1.0.7's small-prime mask and v1.0.8's rewrite ($-42\%$ and $-45\%$ gather thread
instructions at the A100 and H100 pins) each remove instructions without reducing DRAM
traffic or shared atomics. At the A100 pin barrier still leads (34\% on v1.0.7, 33\% on
v1.0.8); at the H100 pin \texttt{long\_scoreboard} leads, 29.4\% against barrier's
19.4\% on v1.0.7 and 25.8\% against 22.1\% on v1.0.8. That is a
property of the pin, not of the architectures: at each device's own tuned tuple
\texttt{long\_scoreboard} leads on both, 32.9\% and 25.4\% against barrier's 18.8 and
15.8\% on v1.0.7, and 36.3\% and 21.0\% against 13.0 and 17.2\% on v1.0.8. A
\texttt{long\_scoreboard} plurality is exposed \emph{latency}, not bandwidth: gather
draws at most 48.4\% of peak DRAM at any of these points.

The profile also exposed a launch-geometry defect invisible to end-to-end timing: the
autotuner drew scatter grids from powers of two, which 108 SMs divide none of --- at
grid 64, 40.7\% of SMs held no work for the kernel taking 34\% of GPU time.
v1.0.6's SM-aligned geometry (scatter grid $64\to108$, gather $512\to864$) lifts
scatter's waves per SM from 0.59 to exactly 1.0 and gather's from 4.74 to 8.0: RSA-100
sieve wall-clock $-17.5\%$ ($-15.1\%$ of total), yield-neutral,
reproduced unprofiled at $n{=}3$ ($-17.4\%$ sieve).

Both fixes transfer only conditionally. SM alignment
\emph{loses} 2.5\% total on the H100: with $2^{k-1}$ polynomials per $a$
(\S\ref{sec:bg}) capping the batch, its RSA-100 configuration ($k=11$) admits no
132-SM-aligned batch above 528, a 48\% cut that halves per-batch amortization, where
the A100 ($k=12$) had headroom to align \emph{upward}.
Attacking gather's barrier bound, v1.0.6's
\emph{halved-footprint} sieving co-resides two gather blocks per SM: a further
$-8.7\%$ total on the A100, $-2.5\%$ on the H100, $+2.2\%$ on the RTX~5070~Ti, at
$+4$--$11\%$ above-idle energy --- it buys GPU-hours, not watt-hours.
 \section{Limitations}
\label{sec:limits}

\subsection{The 2-cycle cliff: genus obstructions}

The large-prime variant has a near-discontinuous failure mode. When the large-prime
relation graph becomes dominated by isolated degree-2 cycles, the genus (class-group)
characters enter the row space and \emph{every} Block Wiedemann dependency becomes
trivial: $X \equiv \pm Y$, so every GCD returns $1$ or $N$. Adding character columns
does \emph{not} fix it --- our packed preprocessor measures the off-diagonal genus
cells at exactly zero. The cliff's location is input- and configuration-dependent: a 94-digit composite succeeds at
large-prime shares up to
${\sim}43\%$ and collapses to a 0\% nontrivial rate at $\ge45.8\%$; RSA-110 is green
across its entire testable range, 51.8--57.1\%, the boundary unlocalized in an
untested 57--62.5\% gap; a known-red configuration exists at 83.6\%; and RSA-155
succeeded at 41.4\% with an 11/24 (45.8\%) nontrivial rate. The widely used
${\sim}55\%$ limit is thus a rule of thumb. To mitigate, keep the share below it, oversieve, and use the matrix-stage
down-filter, which drops large-prime-combined relations above a chosen bound after the
fact, reproducing the relation set a lower-$L$ sieve would have produced \emph{without
re-sieving}.

Relatedly, below ${\sim}85$ digits the single-large-prime variant causes a 100\%
square-root failure rate --- the same obstruction in acute form, the $a$-coefficient
primes and the sieve primes there coming from overlapping regions of the factor base
and creating systematic \GFtwo{} cancellation in combined relations. Such inputs should run with $L=0$; the only automatic gate is a cold-start fallback
below 200 bits, so the intervening range relies on the operator's judgement.

\subsection{Autotuning does not pay for itself within a run}
\label{sec:autotune}

Two optimizations travel under the name ``autotuning''. The cheap part tunes kernel
launch geometry and applies best-known settings from a persistent parameter history.
The expensive part, a joint $(F,L)$ \emph{parameter} search whose probes are real
truncated sieve runs, delivers $1.59\times$ at RSA-100 (305\,s $\to$ 192\,s) but costs
${\sim}300$\,s itself, so it pays only as amortized parameter management across runs.
It is also structurally blind to large primes: probes run at 12\% of target relations
and matching is quadratic in witness count, so a probe observes $0.12^2 = 1.4\%$ of the
full run's matches, the compensating projection model being known-flawed. Probe budgets
overrun (${\sim}80$--160 against a design target of ${\sim}40$), the
linear-algebra cost model was calibrated on one GPU, and the optimizer cannot see the
wide-accumulator path (\S\ref{sec:accum}) used at RSA-150/155 scale --- which is why RSA-155 ran with
autotuning off, whatever it is worth there being upside we do not claim.

\subsection{The 512-bit ceiling, and threats to validity}

The implementation is built on a fixed \texttt{uint512} type, so ${\sim}155$ digits is
a hard architectural limit; RSA-155 also used ${\approx}79\%$ of a 94\,GB H100's
memory, so the largest instance this code can address is close to the largest that
fits --- RSA-155 is its ceiling in both senses. The kernel-level profile of
\S\ref{sec:roofline} covers the
narrow (uint8) accumulator path at RSA-100 on an A100 and an H100; the wide path
carrying RSA-150/155 has no per-kernel profile, so no kernel-level claim extends to
those sizes, and geometry alignment alone should give at most low-single-digit gains
there. RSA-140 is not like-for-like (deduplication safety
factor 1.4 vs 1.05) and is the largest residual in the seven-point fit of \S\ref{sec:benchfact},
and normalization to H100-GPU-hours assumes relative hardware speed at RSA-100
carries over to larger $N$.
 \section{Conclusion and Future Work}
\label{sec:conclusion}

An irregular, branch-heavy factorization pipeline can be made GPU-resident end to end,
competes with the best CPU implementations at RSA-100 and scales to RSA-155 on 64
GPUs; QS linear algebra is effectively solved on GPUs, and sieve throughput is the
frontier. Three extensions are quantified by our own results. Cluster mode is not yet
fully GPU-resident: multi-node large-prime matching runs on the coordinator's CPU
by design (\S\ref{sec:cluster}), and sharding the witness table by large-prime residue would extend the residency claim to
the multi-node case. Every QS size record since 2002 was set with a multiple-large-prime variant
(\S\ref{sec:qsrecords}), whose cycle-finding structure is exactly the irregular,
hash-join-shaped work this paper argues can be made GPU-resident, so implementing it~\cite{lenstra1994twolp,leyland2002threelp} is the
obvious next step. We do not rank it first: the number field
sieve~\cite{lenstra1993,buhler1993} overtakes the quadratic sieve near 100--110 digits
regardless of variant and the gap widens super-polynomially
(\S\ref{sec:cpubaseline} offers one same-machine data point; on GPUs, concurrent work
reports a 155-digit GNFS in about 1.2 GPU-days~\cite{lu2026rsa260}), so multiple large primes
would pay exactly where the quadratic sieve is the wrong algorithm, and nothing here
bears on deployed parameters. Thom\'e's subquadratic lingen~\cite{thome2002lingen} is measured headroom (4\,h\,14\,m
at RSA-155, \S\ref{sec:linalg}) that at 1.6\% of cost does not limit the pipeline; the
per-kernel roofline (\S\ref{sec:roofline}) explains \emph{why} throughput tracks
compute, not bandwidth, and its first diagnosis, SM-aligned launch geometry, is worth $-15\%$ of RSA-100
wall-clock on the A100 and loses $2.5\%$ on the H100 (\S\ref{sec:roofline}).

Two parts carry over to GNFS: the special-q lattice sieve that dominates its cost, whose
inner loops are the same irregular, bucketed scatter/gather pattern made GPU-resident
here~\cite{pollard1993,franke2005,aoki2004bucket}, and its linear algebra, the same
\GFtwo{} Block Wiedemann problem; concurrent work has since
run both on GPUs at record scale~\cite{lu2026rsa260}, but, to our knowledge, no GPU
GNFS implementation is publicly available. A state-of-the-art \emph{open-source} GPU
special-q lattice siever --- the core of a GPU-resident GNFS that carries this paper's
residency and zero-synchronization results over to the sieve, filtering, polynomial
selection~\cite{kleinjung2006} and the algebraic square
root~\cite{couveignes1993,nguyen1998} --- is the best-motivated extension of this work.
GPU sievers are unlikely to be near their ceiling: v1.0.8's binary alone cuts RSA-100 wall-clock by 19--21\% over v1.0.7 on two devices, and 18\% and
36\% on two more across runs (Table~\ref{tab:rsa100}); the lattice
siever's richer structure may hide more.
 
\section*{Acknowledgments}
We gratefully acknowledge computing time provided by the Paderborn Center for
Parallel Computing (PC${}^2$) on the HPC systems Noctua~2 and Otus under the
project ``A multi GPU implementation of the Multiple Polynomial Quadratic
Sieve''. \AIStatement

\section*{Artifact and Data Availability}
\swname{} is open source under LGPL-3.0-only (with a GPLv3 \S7 permission to link the
NVIDIA CUDA Toolkit); the \GFtwo{} Block Wiedemann engine is released separately. Every job,
arm and operating point is identified in the note to Table~\ref{tab:roofline} or in the
artifact's run logs. Table~\ref{tab:benchfact} holds
v1.0.3b--v1.0.5 measurements, Table~\ref{tab:rsa100} v1.0.5 and v1.0.8, the GPU row of Table~\ref{tab:cpubase} v1.0.8;
\S\ref{sec:roofline} profiles v1.0.6--v1.0.8 at the v1.0.6 pins and tuned tuples. Launch-geometry features are flag-gated and default-off
in every release.
  Archived, versioned software releases~\cite{cudampqsall}: v1.0.5~\cite{cudampqs},
  v1.0.6~\cite{cudampqsv106}, v1.0.7~\cite{cudampqsv107},
  v1.0.8~\cite{cudampqsv108} (tag \texttt{99c6196}, source-identical to the measured builds, \texttt{f0868b1} or successors), and
  the \GFtwo{} Block Wiedemann engine~\cite{blockwiedemann}.
  Run data (logs, power samples, profiler exports, manifest): \cite{cudampqsdata}.

\bibliographystyle{IEEEtran}
\bibliography{references}

\begin{thebibliography}{10}
\providecommand{\url}[1]{#1}
\csname url@samestyle\endcsname
\providecommand{\newblock}{\relax}
\providecommand{\bibinfo}[2]{#2}
\providecommand{\BIBentrySTDinterwordspacing}{\spaceskip=0pt\relax}
\providecommand{\BIBentryALTinterwordstretchfactor}{4}
\providecommand{\BIBentryALTinterwordspacing}{\spaceskip=\fontdimen2\font plus
\BIBentryALTinterwordstretchfactor\fontdimen3\font minus
  \fontdimen4\font\relax}
\providecommand{\BIBforeignlanguage}[2]{{%
\expandafter\ifx\csname l@#1\endcsname\relax
\typeout{** WARNING: IEEEtran.bst: No hyphenation pattern has been}%
\typeout{** loaded for the language `#1'. Using the pattern for}%
\typeout{** the default language instead.}%
\else
\language=\csname l@#1\endcsname
\fi
#2}}
\providecommand{\BIBdecl}{\relax}
\BIBdecl

\bibitem{pomerance1985qs}
C.~Pomerance, ``The quadratic sieve factoring algorithm,'' in \emph{Advances in
  Cryptology: Proceedings of EUROCRYPT 84}, ser. Lecture Notes in Computer
  Science, vol. 209.\hskip 1em plus 0.5em minus 0.4em\relax Springer, 1985, pp.
  169--182.

\bibitem{pomerance1996tale}
------, ``A tale of two sieves,'' \emph{Notices of the American Mathematical
  Society}, vol.~43, no.~12, pp. 1473--1485, 1996.

\bibitem{wang2011gpuqs}
\BIBentryALTinterwordspacing
B.~Wang, ``Quadratic sieve on {GPUs},'' MIT, 18.337 course project report,
  Tech. Rep., 2011. [Online]. Available:
  \url{https://courses.csail.mit.edu/18.337/2011/projects/Wang_report.pdf}
\BIBentrySTDinterwordspacing

\bibitem{bernstein2009ecm}
D.~J. Bernstein, T.-R. Chen, C.-M. Cheng, T.~Lange, and B.-Y. Yang, ``{ECM} on
  graphics cards,'' in \emph{Advances in Cryptology --- EUROCRYPT 2009}, ser.
  Lecture Notes in Computer Science, vol. 5479.\hskip 1em plus 0.5em minus
  0.4em\relax Springer, 2009, pp. 483--501.

\bibitem{miele2014cofactorization}
A.~Miele, J.~W. Bos, T.~Kleinjung, and A.~K. Lenstra, ``Cofactorization on
  graphics processing units,'' in \emph{Cryptographic Hardware and Embedded
  Systems --- CHES 2014}, ser. Lecture Notes in Computer Science, vol.
  8731.\hskip 1em plus 0.5em minus 0.4em\relax Springer, 2014, pp. 335--352.

\bibitem{msieve}
J.~Papadopoulos, ``msieve: A library for factoring large integers,''
  \url{https://sourceforge.net/projects/msieve/}, mainline GPU acceleration
  covers NFS polynomial selection stage 1; the quadratic sieve is CPU-only. For
  the CUDA block-Lanczos development line see the NFS@Home fork.

\bibitem{schmidt2011spmv}
B.~Schmidt, H.~Aribowo, and H.-V. Dang, ``Iterative sparse matrix-vector
  multiplication for integer factorization on {GPUs},'' in \emph{Euro-Par 2011
  Parallel Processing}, ser. Lecture Notes in Computer Science, vol.
  6853.\hskip 1em plus 0.5em minus 0.4em\relax Springer, 2011, pp. 413--424.

\bibitem{schmidt2013bwgpu}
------, ``Iterative sparse matrix-vector multiplication for accelerating the
  block {Wiedemann} algorithm over {GF}(2) on multi-graphics processing unit
  systems,'' \emph{Concurrency and Computation: Practice and Experience},
  vol.~25, no.~4, pp. 586--603, 2013.

\bibitem{msievelacuda}
G.~Childers and J.~Papadopoulos, ``msieve\_nfsathome: {NFS@Home} fork of msieve
  with {CUDA} block {Lanczos},''
  \url{https://github.com/gchilders/msieve_nfsathome}, branch
  \texttt{msieve-lacuda-nfsathome}, {CSR} SpMV over GF(2)
  (\texttt{cub/spmv\_engine.cu}): merge-path nonzero partitioning, word-packed
  VBITS-bit vector blocks, \texttt{atomicXor} accumulation, launch parameters
  selected from row-length statistics; device kernels also cover the block
  vector operations (\texttt{common/lanczos/gpu/}). The GPU Lanczos framework
  carries Papadopoulos's msieve public-domain notice; SpMV-engine commits are
  by Childers, 2021--2025. Accessed 2026-08-14.

\bibitem{cavallar2000rsa155}
S.~Cavallar, B.~Dodson, A.~K. Lenstra, W.~Lioen, P.~L. Montgomery, B.~Murphy,
  H.~te~Riele \emph{et~al.}, ``Factorization of a 512-bit {RSA} modulus,'' in
  \emph{Advances in Cryptology --- EUROCRYPT 2000}, ser. Lecture Notes in
  Computer Science, vol. 1807.\hskip 1em plus 0.5em minus 0.4em\relax Springer,
  2000, pp. 1--18.

\bibitem{weis2026rsa896}
\BIBentryALTinterwordspacing
S.~A. Weis, ``{RSA-896},'' Blog post, Sep. 2026, factored 2026-09-19 (270
  digits, 896 bits) with a GPU port of CADO-NFS made with Claude; up to 2{,}048
  GPUs, about 30 GPU-years over 10 days; stage split not disclosed. Accessed
  2026-09-22. [Online]. Available: \url{https://saweis.net/posts/rsa-896.html}
\BIBentrySTDinterwordspacing

\bibitem{silverman1987mpqs}
R.~D. Silverman, ``The multiple polynomial quadratic sieve,'' \emph{Mathematics
  of Computation}, vol.~48, no. 177, pp. 329--339, 1987.

\bibitem{contini1997siqs}
S.~P. Contini, ``Factoring integers with the self-initializing quadratic
  sieve,'' Master's thesis, University of Georgia, 1997.

\bibitem{crandall2005primes}
R.~Crandall and C.~Pomerance, \emph{Prime Numbers: A Computational
  Perspective}, 2nd~ed.\hskip 1em plus 0.5em minus 0.4em\relax Springer, 2005.

\bibitem{heinrichs2026}
C.~Heinrichs, ``An implementation of the quadratic sieve in {CUDA},'' Master's
  thesis, Paderborn University, Feb. 2026.

\bibitem{dixon1981}
J.~D. Dixon, ``Asymptotically fast factorization of integers,''
  \emph{Mathematics of Computation}, vol.~36, no. 153, pp. 255--260, 1981.

\bibitem{davis1984factorization}
J.~A. Davis and D.~B. Holdridge, ``Factorization using the quadratic sieve
  algorithm,'' in \emph{Advances in Cryptology: Proceedings of Crypto
  83}.\hskip 1em plus 0.5em minus 0.4em\relax Plenum Press, 1984, pp. 103--113.

\bibitem{pomerance1988pipeline}
C.~Pomerance, J.~W. Smith, and R.~Tuler, ``A pipeline architecture for
  factoring large integers with the quadratic sieve algorithm,'' \emph{SIAM
  Journal on Computing}, vol.~17, no.~2, pp. 387--403, 1988.

\bibitem{alford1995siqs}
W.~R. Alford and C.~Pomerance, ``Implementing the self-initializing quadratic
  sieve on a distributed network,'' in \emph{Number-Theoretic and Algebraic
  Methods in Computer Science (Moscow, 1993)}.\hskip 1em plus 0.5em minus
  0.4em\relax World Scientific, 1995, pp. 163--174.

\bibitem{peralta1995hypercube}
R.~Peralta, ``Implementation of the hypercube multiple polynomial quadratic
  sieve,'' 1995.

\bibitem{knuth2011taocp4a}
D.~E. Knuth, \emph{The Art of Computer Programming, Volume 4A: Combinatorial
  Algorithms, Part 1}.\hskip 1em plus 0.5em minus 0.4em\relax Addison-Wesley,
  2011.

\bibitem{lenstra1994twolp}
A.~K. Lenstra and M.~S. Manasse, ``Factoring with two large primes,''
  \emph{Mathematics of Computation}, vol.~63, no. 208, pp. 785--798, 1994.

\bibitem{leyland2002threelp}
P.~Leyland, A.~K. Lenstra, B.~Dodson, A.~Muffett, and S.~S. {Wagstaff Jr.},
  ``{MPQS} with three large primes,'' in \emph{Algorithmic Number Theory
  (ANTS-V)}, ser. Lecture Notes in Computer Science, vol. 2369.\hskip 1em plus
  0.5em minus 0.4em\relax Springer, 2002, pp. 446--460.

\bibitem{coppersmith1994bw}
D.~Coppersmith, ``Solving homogeneous linear equations over {GF}(2) via block
  {Wiedemann} algorithm,'' \emph{Mathematics of Computation}, vol.~62, no. 205,
  pp. 333--350, 1994.

\bibitem{wiedemann1986}
D.~H. Wiedemann, ``Solving sparse linear equations over finite fields,''
  \emph{IEEE Transactions on Information Theory}, vol.~32, no.~1, pp. 54--62,
  1986.

\bibitem{kaltofen1995bw}
E.~Kaltofen, ``Analysis of {Coppersmith}'s block {Wiedemann} algorithm for the
  parallel solution of sparse linear systems,'' \emph{Mathematics of
  Computation}, vol.~64, no. 210, pp. 777--806, 1995.

\bibitem{villard1997bw}
G.~Villard, ``Further analysis of {Coppersmith}'s block {Wiedemann} algorithm
  for the solution of sparse linear systems,'' in \emph{Proceedings of ISSAC
  '97}.\hskip 1em plus 0.5em minus 0.4em\relax ACM, 1997, pp. 32--39.

\bibitem{montgomery1995}
P.~L. Montgomery, ``A block {L}anczos algorithm for finding dependencies over
  {GF(2)},'' in \emph{Advances in Cryptology --- EUROCRYPT '95}, ser. Lecture
  Notes in Computer Science, vol. 921.\hskip 1em plus 0.5em minus 0.4em\relax
  Springer, 1995, pp. 106--120.

\bibitem{geiselmann2006}
W.~Geiselmann, F.~Januszewski, H.~K{\"o}pfer, J.~Pelzl, and R.~Steinwandt, ``A
  simpler sieving device: Combining {ECM} and {TWIRL},'' in \emph{Information
  Security and Cryptology --- ICISC 2006}, ser. Lecture Notes in Computer
  Science, vol. 4296.\hskip 1em plus 0.5em minus 0.4em\relax Springer, 2006,
  pp. 118--135.

\bibitem{cadonfs}
{The CADO-NFS Development Team}, ``{CADO-NFS}: An implementation of the number
  field sieve algorithm,'' \url{https://cado-nfs.gitlabpages.inria.fr/},
  repository \url{https://gitlab.inria.fr/cado-nfs/cado-nfs}, the development
  branch additionally contains a CPU SIQS by C.~Bouvier:
  \texttt{siqs}/\texttt{siqs\_tracek} binaries building the \texttt{las} siever
  source with \texttt{SIQS\_SIEVE=1} (commit \texttt{7486519}, 2025-12-18);
  Gray-code enumeration and level-2 bucket sieving (January 2026); SSE2/AVX2
  kernels and \texttt{cadofactor} integration with parameters for 40--90-digit
  composites (merge request \texttt{cado-nfs!260}, merged 2026-07-06). Not in
  any tagged release. Accessed 2026-08-14.

\bibitem{lu2026rsa260}
\BIBentryALTinterwordspacing
E.~Lu, ``Factoring {RSA-260},'' Cognition blog, Sep. 2026, published
  2026-09-09; factors found 2026-09-03. GNFS with a modified CADO-NFS developed
  with the Devin coding agent; GPU polynomial selection (stage~1), lattice
  sieving, Block Wiedemann and square root, filtering on CPU; about 4{,}900
  GPU-days (NVIDIA B200/GB200/GB300). Accessed 2026-09-22. [Online]. Available:
  \url{https://cognition.com/blog/factoring-rsa-260}
\BIBentrySTDinterwordspacing

\bibitem{atkins1995squeamish}
D.~Atkins, M.~Graff, A.~K. Lenstra, and P.~C. Leyland, ``The magic words are
  squeamish ossifrage,'' in \emph{Advances in Cryptology --- ASIACRYPT '94},
  ser. Lecture Notes in Computer Science, vol. 917.\hskip 1em plus 0.5em minus
  0.4em\relax Springer, 1995, pp. 261--277.

\bibitem{konsor2020rsa140qs}
\BIBentryALTinterwordspacing
P.~Konsor, ``Useless accomplishment: {RSA-140} factorization with quadratic
  sieve,'' mersenneforum.org thread 25676, 2020, {QS} record: 140 digits,
  ${\sim}6{,}000$ core-hours. [Online]. Available:
  \url{https://www.mersenneforum.org/showthread.php?t=25676}
\BIBentrySTDinterwordspacing

\bibitem{buhrow2025rsa150qs}
\BIBentryALTinterwordspacing
B.~Buhrow, ``Three large primes (again) and another useless accomplishment,''
  mersenneforum.org, post 1079832, Jun. 2025, posted 2025-06-18 as user
  \texttt{bsquared}. YAFU three-large-prime (TLP) variant on AMD EPYC 9174
  processors; cost in core-hours: RSA-100 0.96, RSA-110 6.7, RSA-120 40.9,
  RSA-130 256, RSA-140 1{,}791, RSA-150 11{,}664 (LPB 32 bits, FB bound 1.44M).
  RSA-150 is the QS record prior to this work; its linear algebra used block
  Lanczos. [Online]. Available:
  \url{https://www.mersenneforum.org/node/1079832#post1079832}
\BIBentrySTDinterwordspacing

\bibitem{yafu}
------, ``{YAFU}: Yet another factorization utility,''
  \url{https://github.com/bbuhrow/yafu}, {AVX-512} {SIQS}; holds prior {QS}
  records including {RSA-150} (June 2025).

\bibitem{kleinjung2010rsa768}
T.~Kleinjung, K.~Aoki, J.~Franke, A.~K. Lenstra, E.~Thom{\'e}, J.~W. Bos,
  P.~Gaudry, A.~Kruppa, P.~L. Montgomery, D.~A. Osvik, H.~te~Riele,
  A.~Timofeev, and P.~Zimmermann, ``Factorization of a 768-bit {RSA} modulus,''
  in \emph{Advances in Cryptology --- CRYPTO 2010}, ser. Lecture Notes in
  Computer Science, vol. 6223.\hskip 1em plus 0.5em minus 0.4em\relax Springer,
  2010, pp. 333--350.

\bibitem{boudot2020rsa250}
F.~Boudot, P.~Gaudry, A.~Guillevic, N.~Heninger, E.~Thom{\'e}, and
  P.~Zimmermann, ``Comparing the difficulty of factorization and discrete
  logarithm: A 240-digit experiment,'' in \emph{Advances in Cryptology ---
  CRYPTO 2020}, ser. Lecture Notes in Computer Science, vol. 12171.\hskip 1em
  plus 0.5em minus 0.4em\relax Springer, 2020, pp. 62--91.

\bibitem{factorrecords}
\BIBentryALTinterwordspacing
P.~Zimmermann, ``Integer factoring records,'' Online database of documented
  GNFS factorizations. [Online]. Available:
  \url{https://members.loria.fr/PZimmermann/records/factor.html}
\BIBentrySTDinterwordspacing

\bibitem{aoki2004bucket}
K.~Aoki and H.~Ueda, ``Sieving using bucket sort,'' in \emph{Advances in
  Cryptology --- ASIACRYPT 2004}, ser. Lecture Notes in Computer Science, vol.
  3329.\hskip 1em plus 0.5em minus 0.4em\relax Springer, 2004, pp. 92--102.

\bibitem{franke2005}
J.~Franke and T.~Kleinjung, ``Continued fractions and lattice sieving,'' in
  \emph{Proc.\ Special-purpose Hardware for Attacking Cryptographic Systems
  (SHARCS)}, 2005.

\bibitem{montgomery1985}
P.~L. Montgomery, ``Modular multiplication without trial division,''
  \emph{Mathematics of Computation}, vol.~44, no. 170, pp. 519--521, 1985.

\bibitem{thome2002lingen}
E.~Thom{\'e}, ``Subquadratic computation of vector generating polynomials and
  improvement of the block {Wiedemann} algorithm,'' \emph{Journal of Symbolic
  Computation}, vol.~33, no.~5, pp. 757--775, 2002.

\bibitem{ding2019roofline}
N.~Ding and S.~Williams, ``An instruction roofline model for {GPUs},'' in
  \emph{2019 IEEE/ACM Performance Modeling, Benchmarking and Simulation of High
  Performance Computer Systems (PMBS)}.\hskip 1em plus 0.5em minus 0.4em\relax
  IEEE, 2019, pp. 7--18.

\bibitem{lenstra1993}
A.~K. Lenstra and J.~H.~W.~Lenstra, Eds., \emph{The Development of the Number
  Field Sieve}, ser. Lecture Notes in Mathematics.\hskip 1em plus 0.5em minus
  0.4em\relax Springer, 1993, vol. 1554.

\bibitem{buhler1993}
J.~P. Buhler, J.~H.~W.~Lenstra, and C.~Pomerance, ``Factoring integers with the
  number field sieve,'' in \emph{The Development of the Number Field Sieve},
  ser. Lecture Notes in Mathematics.\hskip 1em plus 0.5em minus 0.4em\relax
  Springer, 1993, vol. 1554, pp. 50--94.

\bibitem{pollard1993}
J.~M. Pollard, ``The lattice sieve,'' in \emph{The Development of the Number
  Field Sieve}, ser. Lecture Notes in Mathematics.\hskip 1em plus 0.5em minus
  0.4em\relax Springer, 1993, vol. 1554, pp. 43--49.

\bibitem{kleinjung2006}
T.~Kleinjung, ``On polynomial selection for the general number field sieve,''
  \emph{Mathematics of Computation}, vol.~75, no. 256, pp. 2037--2047, 2006.

\bibitem{couveignes1993}
J.-M. Couveignes, ``Computing a square root for the number field sieve,'' in
  \emph{The Development of the Number Field Sieve}, ser. Lecture Notes in
  Mathematics.\hskip 1em plus 0.5em minus 0.4em\relax Springer, 1993, vol.
  1554, pp. 95--102.

\bibitem{nguyen1998}
P.~Nguyen, ``A {M}ontgomery-like square root for the number field sieve,'' in
  \emph{Algorithmic Number Theory (ANTS-III)}, ser. Lecture Notes in Computer
  Science, vol. 1423.\hskip 1em plus 0.5em minus 0.4em\relax Springer, 1998,
  pp. 151--168.

\bibitem{cudampqsall}
\BIBentryALTinterwordspacing
F.~Januszewski and C.~Heinrichs, ``{CUDA-MPQS}: A fully {GPU}-resident
  quadratic sieve,'' Computer software, all versions (concept DOI), Zenodo,
  2026, source at \url{https://github.com/drjanosch42/cuda-mpqs}. [Online].
  Available: \url{https://doi.org/10.5281/zenodo.21619468}
\BIBentrySTDinterwordspacing

\bibitem{cudampqs}
\BIBentryALTinterwordspacing
------, ``{CUDA-MPQS}: A fully {GPU}-resident quadratic sieve,'' Computer
  software, version v1.0.5, Zenodo, 2026, open source (LGPL-3.0-only with
  NVIDIA CUDA Toolkit exception); results in this paper correspond to v1.0.5f;
  source at \url{https://github.com/drjanosch42/cuda-mpqs}. [Online].
  Available: \url{https://doi.org/10.5281/zenodo.21619473}
\BIBentrySTDinterwordspacing

\bibitem{cudampqsv106}
\BIBentryALTinterwordspacing
------, ``{CUDA-MPQS}: A fully {GPU}-resident quadratic sieve,'' Computer
  software, version v1.0.6, Zenodo, 2026, open source (LGPL-3.0-only with
  NVIDIA CUDA Toolkit exception); the bottleneck analysis
  (\S\ref{sec:roofline}) profiles and launch-geometry comparison use this
  release; source at \url{https://github.com/drjanosch42/cuda-mpqs}. [Online].
  Available: \url{https://doi.org/10.5281/zenodo.22165743}
\BIBentrySTDinterwordspacing

\bibitem{cudampqsv107}
\BIBentryALTinterwordspacing
------, ``{CUDA-MPQS}: A fully {GPU}-resident quadratic sieve,'' Computer
  software, version v1.0.7, Zenodo, 2026, open source (LGPL-3.0-only with
  NVIDIA CUDA Toolkit exception); the v1.0.7 roofline columns and the v1.0.7
  baselines of the v1.0.8 comparisons use this release; source at
  \url{https://github.com/drjanosch42/cuda-mpqs}. [Online]. Available:
  \url{https://doi.org/10.5281/zenodo.22902652}
\BIBentrySTDinterwordspacing

\bibitem{cudampqsv108}
\BIBentryALTinterwordspacing
------, ``{CUDA-MPQS}: A fully {GPU}-resident quadratic sieve,'' Computer
  software, version v1.0.8, Zenodo, 2026, open source (LGPL-3.0-only with
  NVIDIA CUDA Toolkit exception); the v1.0.8 RSA-100 column of
  Table~\ref{tab:rsa100}, the GPU row of Table~\ref{tab:cpubase} and the v1.0.8
  roofline column use this release (measured on \texttt{f0868b1} or
  source-identical builds; the GitHub tag v1.0.8, commit \texttt{99c6196},
  differs only in documentation and metadata); source at
  \url{https://github.com/drjanosch42/cuda-mpqs}. [Online]. Available:
  \url{https://doi.org/10.5281/zenodo.23137108}
\BIBentrySTDinterwordspacing

\bibitem{blockwiedemann}
\BIBentryALTinterwordspacing
F.~Januszewski, ``A {GPU}-resident block {Wiedemann} {GF(2)} linear-algebra
  implementation,'' Computer software, version v1.0.2, Zenodo, 2026, open
  source (LGPL-3.0-only); Stage-2 initialization on device since v1.0.2
  (bundled with \swname{} v1.0.8; v1.0.5 bundles v1.0.1); source at
  \url{https://github.com/drjanosch42/block-wiedemann}. [Online]. Available:
  \url{https://doi.org/10.5281/zenodo.23136111}
\BIBentrySTDinterwordspacing

\bibitem{cudampqsdata}
\BIBentryALTinterwordspacing
F.~Januszewski and C.~Heinrichs, ``{CUDA-MPQS}: Run data for the {RSA-100} to
  {RSA-155} benchmark factorizations,'' Dataset, version v1.0, CC BY 4.0,
  Zenodo, 2026, raw run logs, power samples, launcher scripts, profiler reports
  and CSV exports, analysis scripts, and a sha256 manifest behind every
  reported figure (\texttt{cuda-mpqs-data-v1.0.tar.gz}). [Online]. Available:
  \url{https://doi.org/10.5281/zenodo.23144292}
\BIBentrySTDinterwordspacing

\end{thebibliography}

\end{document}